\documentclass[12pt]{amsart}
\usepackage{amsmath, amscd, amssymb, latexsym}

\theoremstyle{plain}
\newtheorem{thm}{Theorem}[section]
\newtheorem{lem}[thm]{Lemma}
\newtheorem{cor}[thm]{Corollary}
\newtheorem{pro}[thm]{Proposition}

\theoremstyle{definition}
\newtheorem{df}[thm]{Definition}

\def\cok{\mbox{coker}}
\def\om{\omega}
\def\Om{\Omega}
\def\La{\Lambda}

\def\ov{\overline}
\def\mod{\mbox{ mod}}
\def\al{\alpha}
\def\ga{\gamma}
\def\Ga{\Gamma}
\def\si{\sigma}
\def\Si{\Sigma}
\def\ep{\epsilon}

\def\wt{\widetilde}

\def\nt{\noindent}

\def\Dr{{\not}\partial}
\def\lra{\longrightarrow}

\def\pa{\partial}

\def\de{\delta}
\def\ti{\tilde}
\def\mc{\mathcal}
\def\mb{\mathbb}
\def\i{{\bf i}}
\def\mr{\mathrm}
\def\R{\mr{R}}

\begin{document}

\hspace{-30mm}

\nt{Intern. J. Math. (2009)}

\vspace{10mm}

\title{Orientability and real Seiberg-Witten invariants}

\author[G. Tian]{Gang Tian}
\address{Department of Mathematics\\
Princeton University\\Princeton\\
NJ 08544} \email{tian@math.princeton.edu}

\author[S. Wang]{Shuguang Wang}
\address{Department of Mathematics\\University of Missouri\\ Columbia\\
MO 65211} \email{wangs@missouri.edu}

\begin{abstract}
 We investigate  Seiberg-Witten theory in the presence of real
structures. Certain conditions are obtained so that integer valued real
Seiberg-Witten invariants can be defined. In general we study
properties of the real Seiberg-Witten projection map from the
point of view of Fredholm map degrees.
  \end{abstract}

\maketitle

\noindent 2000 MS Classification: 57R57, 53C07, 14P25

\section{Introduction}\label{intro}

After the much success of  Seiberg-Witten theory, it is a natural
problem to study the real version and its potential application in
real algebraic geometry. In outline, the real version starts with
a K\"{a}hler 4-manifold with anti-holomorphic involution (a real
structure). Then one would like to understand the lifted action on
the Seiberg-Witten moduli space as well as the invariant extracted
from the real moduli space. This real Seiberg-Witten invariant
should link to the real Gromov-Witten invariant counting real
holomorphic curves through a Taubes type correspondence as in
\cite{t2}. An invariant counting nodal rational curves in real
rational surfaces is  found in \cite{we}. More recently Solomon
\cite{so} has defined  Gromov-Witten invariants counting arbitrary
degree real holomorphic curves from fixed Riemann surfaces. A
K\"{a}hler manifold with a real structure is what physicists refer
to as an orientifold \cite{sv, bfm}, and the real Gromov-Witten
invariants have been one of their main interests for the last few
years. Compare with the original Gromov-Witten theory of Ruan-Tian
\cite{rt}.

One should first point  out that such a real theory is not to be
treated as an equivariant theory; for one thing, the lifted real map
does not preserve the spin$^c$ bundle in the Seiberg-Witten
theory. Nevertheless the lifted map is conjugate linear in a
proper sense. Consequently, among the standard issues of
transversality, compactness, orientability and reducible
solutions, only orientability requires a substantially  new
strategy to tackle. As a matter of fact, the real Seiberg-Witten
moduli space  is not necessarily orientable or naturally oriented
even if orientable. This is rather typical in real algebraic
geometry: the real part of a real structure is usually
non-orientable or un-oriented.

It is well-known that the usual Seiberg-Witten invariant can be
viewed as the degree of the projection map $\pi:{\bf M}\to
\i\Om^2_+$, where ${\bf M}$ is the parameterized moduli space.
This is the case if the moduli space is 0-dimensional. In the real
Seiberg-Witten theory, we will encounter the real Seiberg-Witten
projection map $\pi_{\R}:{\bf M}_{\R}\to (\i \Om^2_+)_{\R}$
defined on the real parameterized moduli space. We will undertake
two approaches: the first is to place  real moduli spaces in the
real configuration space $\mc{B}_\R$ and seek conditions in terms
of $\mc{B}_\R$ that will guarantee the orientability of real
moduli spaces. To this end we have the following results (see
Theorems \ref{prsw} and \ref{isw}).

\begin{thm}
Let $X$ be a K\"{a}hler surface with a real structure $\si$
and $S$ a spin$^c$ structure compatible with $\si$.
Fix orientations on  $H^1_\R(X,{\bf R}), H^+_\R(X,{\bf R})$.
If $H^1(X,{\bf R})$ is trivial or if $c_1(L)$ is divisible by 4 for
the determinant $L$ of $S$,
then the real Seiberg-Witten invariant is  well-defined and takes
integer values.
\end{thm}

As an application we prove a real version of the Thom conjecture
(Corollary \ref{tho}) for smoothly embedded surfaces in ${\bf
CP}^2$ that are equivariant with respect to the real structures.

The second approach is less conventional,  where we focus on the
parameterized moduli space ${\bf M}_{\R}$ itself without ever
involving  $\mc{B}_\R$. Though as a trade-off, we need to work
with all perturbations in $(\i \Om^2_+)_{\R}$.
 The main goal here is to understand the critical point set and
the regular value set of the projection $\pi_{\R}$. Since
$\pi_{\R}$ is proper, its regular values form an open and dense
subset of $(\i \Om^2_+)_{\R}$. Thus the complement forms
``walls'', cutting the regular value set into chambers.  In the
absence of the orientability and hence the integer Seiberg-Witten
invariant, the pattern of chambers and distribution of the
chamber-wide Seiberg-Witten invariants become new geometry to
investigate for our real Seiberg-Witten theory. Among the main
results here, we prove the following (Theorems \ref{cr} and
\ref{crr})

\begin{thm}
Let ${\bf C_\R}$
 denote the critical point set of $\pi_\R$ and
 ${\bf C}_\R(l)=\{{\bf x}\in{\bf C_\R}\mid 
\dim\mr{coker}D\pi_\R({\bf x})=l\}$.
For each integer $l\geq 0$, ${\bf C}_\R(l)\subset {\bf M}_\R$ is a
smooth Banach submanifold of codimension $kl$, 
where $k=\mr{ind}D\pi_\R+l$.
 \end{thm}

Here is an outline of the paper. In Section \ref{per}, after reviewing
the set-up and notations of the standard Seiberg-Witten projection,
  we discuss thoroughly how to lift a real structure from
an almost complex manifold to its associated spin$^c$ bundle, and
apply the lifted real structure to the Seiberg-Witten theory. In
Section \ref{dlin}, we determine the orientation bundle of real moduli
spaces and examine the natural extension to the real configuration
space.
 This illustrates precisely the difference between the usual and real
Seiberg-Witten theories. In Section \ref{rss} we find sufficient
conditions so that the real moduli spaces are orientable and oriented,
thus defining integer valued real invariants. In Section \ref{yup},
we demonstrate that the
critical point set of $\pi_{\R}$ stratifies into immersed
submanifolds of ${\bf M}_{\R}$, which are of the expected
co-dimensions.  Under the assumption that $\pi_{\R}$ is
non-orientable, we introduce chamberwise invariants and
their distribution.

\vspace{5mm}

\nt{\bf Acknowledgments.} The project was initiated in 2004 when
both authors were visiting MSRI. We thank the institute for
providing the excellent environment. Work was partially supported
by NSF grants.

\vspace{2mm}

\section{Lifted real structures and Seiberg-Witten 
equations}\label{per}

\subsection{Parameterized Seiberg-Witten moduli spaces}\label{prel}
 We first recall briefly the standard Seiberg-Witten theory
and set up notations to be used; compare for example
  \cite{m, sa}. Special care is placed on the differential of the
projection map into the perturbation space. The calculations are
often left out in the literature, partly because of the similarity
with the previous Donaldson theory.

Start with the perturbed Seiberg-Witten equations for a
general spin$^c$ structure $S=S^+\oplus S^-$ with determinant $L$,
defined on an arbitrary smooth 4-manifold $X$. The equations are:
 \begin{equation}\label{SW}
 \begin{array}{l}
\Dr_A\Phi=0\\
 F^+_A=q(\Phi)-h
\end{array}
\end{equation}
 where $A$ is a connection on $L$, $\Phi\in\Ga (S^+)$,
and $h\in\i\Omega^2_+$ is
 a perturbation.

From now on, we will suppress Sobolev spaces throughout the paper in order
to focus on the main issue of orientability.
 Given any solution
$(A,\Phi)$ of (\ref{SW}), we have the following fundamental elliptic complex
\begin{equation}\label{complex}
 \begin{array}{ccccc}
&D^1&&D^2&\\
0\lra\i \Omega^0&\lra&\i\Omega^1\oplus\Ga(S^+)&\lra&
\i\Omega^2_+\oplus\Ga(S^-)\lra0
\end{array}\end{equation}
where $D^1=D^1_{A,\Phi}$ is defined by $D^1(f)=(2df, -f\Phi)$ for
$f\in \i\Omega^0$, and $D^2=D^2_{A,\Phi}$ by
$$D^2(a,\phi)=(d^+a-Dq_{\Phi}(\phi),
\Dr_A\phi+\frac{1}{2}a\cdot\Phi)$$ for
$(a,\phi)\in\i\Omega^1\oplus\Ga(S^+)$. Here $q(\Phi)$ has the
differential
 $$Dq_\Phi(\phi)=\Phi\otimes\phi^*+\phi\otimes\Phi^*-
 \frac{<\Phi,\phi>+\ov{<\Phi,\phi>}}{2}\text{Id}.$$
 Remark that only the first equation $\Dr_A\Phi=0$ is needed to
show that (\ref{complex}) is a complex. Moreover
 the perturbation $h$ does
not appear explicitly in the formulas but certainly affects $D^1, D^2$ through
$A, \Phi$ due to equations (\ref{SW}).

The perturbed SW equations (\ref{SW}) define the smooth function
\begin{equation}\label{swf}
\begin{array}{cccc}
\mc{F}:&\mc{A}\times\Ga(S^+)\times\i\Om^2_+&\lra&\i\Om^2_+\times\Ga(S^-)\\
&(A,\Phi, h)&\mapsto&(F^+_A-q(\Phi)-h, \Dr_A\Phi).
\end{array}\end{equation}
At a point $(A,\Phi,h)$, the differential $D\mc{F}_{A,\Phi,h}:\;
\i\Om^1\oplus\Ga(S^+)\oplus\i\Om^2_+ \lra\i\Om^2_+\oplus\Ga(S^-)$
 is
 \begin{equation}\label{DF}
  D\mc{F}(a,\phi,k)=D^2(a,\phi)-(k,0).
  \end{equation}
   Here $(k,0)$ is viewed
as a vector in the direct sum. The standard transversality theorem says that
$0$ is
a regular value of $\mc{F}$, when restricted to irreducibles
$(A,\Phi,h), \Phi\not=0$. Hence $\mc{F}^{-1}(0)\cap\{\Phi\not=0\}$ is
a smooth
Banach manifold. The tangent space at such a point $(A,\Phi,h)$ is of curse
$T\mc{F}^{-1}(0)=\ker D\mc{F}_{A,\Phi,h}$.

Take the projection to the parameter space $\pi:
\mc{F}^{-1}(0)\to\i\Om^2_+$, namely $\pi(A,\Phi,h)=h$. We want to express
the kernel and cokernel of its differential $D\pi$ in terms
of $D^2$ introduced
above. The differential
at a point $(A,\Phi,h)$ with $\Phi\not=0$ is
$$\begin{array}{cccc}
D\pi:&\ker D\mc{F}_{A,\Phi,h}&\lra&\i\Om^2_+\\
&(a,\phi,k)&\mapsto&k.\end{array}$$
Then one can
readily verify there is a natural isomorphism $\ker D\pi=\ker D^2$, by
using (\ref{DF}) and the inclusion
 $\i\Om^1\oplus\Ga(S^+)\hookrightarrow\i\Om^1\oplus\Ga(S^+)\oplus\i\Om^2_+$.
  Next one  can relate the image sets {im}$D\pi\subset\i\Om^2_+$ and
 im$D^2\subset \i\Omega^2_+\oplus\Ga(S^-)$ as follows:
 \begin{equation}
 \text{im}D\pi=\{p_1(\xi)\mid \xi\in\text{im} D^2 \text{ such that }
  p_2(\xi)=0\},
 \end{equation}
where $p_1,p_2$ are projections of $\i\Om^2_+\oplus \Ga(S^-)$ onto
its factors. Thus the inclusion map $\i\Om^2_+\hookrightarrow
\i\Om^2_+\oplus\Ga(S^-), k\mapsto(k,0)$ induces a well-defined
injective map
$$\i\Om^2_+/\text{im}D\pi\lra
\i\Om^2_+\oplus\Ga(S^-)/\text{im}D^2.$$
 Furthermore this map is surjective, which follows from
 (\ref{DF}) and
  $\text{coker} D\mc{F}=0$ on $\mc{F}^{-1}(0)$ by the transversality
theorem.
Hence we have a natural isomorphism
\begin{equation}
 \text{coker}D\pi_{A,\Phi,h}\lra \text{coker}D^2_{A,\Phi},
 \end{equation}
 induced by the inclusion $\i\Om^2_+\hookrightarrow
 \i\Om^2_+\oplus \Ga(S^-)$.
Again the perturbation $h$ does  affect $\text{coker}D^2_{A,\Phi}$
through the SW equations.

Consider the gauge $\mc{G}$ action, $g(A, \Phi)=((g^2)^*A,
g\cdot\Phi)$ for each $g\in\mc{G}$. After modulo out the action,
we have the parameterized moduli space ${\bf  M}
=\mc{F}^{-1}(0)/\mc{G}$, which is however only a topological
space. But the subspace ${\bf  M}^*$ of irreducible solutions is a
smooth Banach manifold and  the tangent space at a point
$[A,\Phi,h]$ is
 $$
 \mr{T}{\bf  M}^* = \ker D\mc{F}\cap\ker (D^1)^*,
 $$
by using the slice $\ker (D^1)^*$ of the $\mc{G}$ action.
(It will be useful to keep in mind that
$(D^1)^*(a,\phi)=2d^*a-<\phi,\Phi>$.) We have the new
projection map $\pi:{\bf M} \lra \i\Om^2_+$. By applying
the slice to other discussions above we can
summarize the main results here.

\begin{pro}\label{dpi}
 (i) The tangent space of the parameterized irreducible moduli
 space ${\bf  M}^* $
 at a point $[A,\Phi,h]$ is the following subspace of
$\i\Om^1\oplus\Ga(S^+)\oplus\i\Om^2_+$:
 $$\mr{T}{\bf  M}^* =\{(a,\phi,k)\mid D^2(a,\phi)=(k,0) \text{ and }
(D^1)^*(a,\phi)=0\}.$$

(ii) The differential $D\pi:\mr{T}{\bf  M} ^*\to\i\Om^2_+$
 has the kernel and cokernel
canonically identified with:
$$\ker(D\pi)={\bf H}^1, \; \mr{coker}D\pi={\bf H}^2,$$
where ${\bf H}^1={\bf H}^1_{A,\Phi}, {\bf H}^2={\bf H}^2_{A,\Phi}$
are the cohomology of the complex (\ref{complex}).
\end{pro}

Because of the slice $\ker(D^1)^*$, we cannot use $\ker D^2$ alone to
characterize $\ker D\pi$ on $\mr{T}{\bf M}^*$, although we do have
$\mr{coker}D^2=\mr{coker}D\pi$.

Quite often, it is advantageous to
form a single elliptic operator converted from the basic complex
(\ref{complex}):
\begin{equation}\label{de} \de=\de_{A,\Phi}=D^2\oplus
(D^1)^*:\i\Om^1\oplus\Ga(S^+)\lra [\i\Om^2_+\oplus
\Ga(S^-)]\oplus\i\Om^0.
 \end{equation}
Remark that for $[A,\Phi, h]\notin{\bf M}^*$, (\ref{complex}) is not a
complex but $\de=\de_{A,\Phi}$ is still elliptic.

 At $[A,\Phi, h]\in{\bf M}^*$, $\ker\de={\bf H}^1$,
 $\mr{coker}\de={\bf H}^2$. Hence
  Proposition \ref{dpi} translates into the following:
\begin{cor}\label{dee} There are natural isomorphisms
 \begin{equation}\label{de0}
  \ker D\pi=
 \ker\de,\; \mr{coker}D\pi=\mr{coker}\de.
\end{equation}
\end{cor}

Unlike $D\pi$, $\de$ is defined between two fixed vector bundles,
i.e. $\de$ can be viewed as a family of elliptic operators, which
is another advantage over $D\pi$.

\subsection{Real spin$^c$ structures}
Let $(X, J)$ be an almost complex manifold of dimension $2n$ and
$\si: X\to X$ a real structure, i.e.  an anti-holomorphic
involution, so $\si_*J=-J\si_*:TX\to TX$. Endow $X$ with a Hermitian metric
that is preserved by both $J$ and $\si$, namely,
$(Ju, Jv)=(u,v), (\si_*(u),\si_*(v))=(u,v)$.
It is
well-known that $X$ has a
canonical spin$^c$ structure $P_{sp}$ associated with $J$ and the metric.
In this subsection, we consider a natural lifting of $\si$ on the
spin$^c$ structure.

Let $P_U\lra X$ be the $U(n)$-bundle of complex frames
and $P_{so}\lra X$ the $SO(2n)$-bundle of real frames, both of which use
the metric on $X$. There is a
natural inclusion $\rho: U(n)\lra SO({2n})$, given as
$\rho[a_{jk}]=[a_{jk, r}]$, where
$$a_{jk,r}=\left[\begin{array}{rr} x&y\\-y&x\end{array}\right]$$
if the entry $a_{jk}=x+{\bf i}y$. Then
$P_{so}=P_U\times_{\rho}SO({2n})$. Note $\rho(\ov{u})=T\rho(u)T$,
where  $T=T^{-1}$ is the diagonal $2n\times 2n$ matrix
\begin{equation}\label{dt}
\text{diag}\{1,-1,1,-1,\cdots,1,-1\}.
\end{equation}
 Hence
 under $\rho$, the complex
conjugation  on $U(n)$ is transferred onto $SO(2n)$ as
$v\mapsto\ov{v}:=TvT$.

\begin{lem}\label{u} There is a canonical involution lifting
$\tau: P_U\lra P_U$ of $\si$ which is conjugate in the sense that
$\tau(pu)=\tau(p)\ov{u}$ for $p\in P_U, u\in U({n})$. Moreover,
$\tau$ induces a lifting on $P_{so}$ (still denoted by $\tau$)
satisfying $\tau(pv)=\tau(p)\ov{v}$ for $p\in P_{so}, v\in
SO(2n)$.
\end{lem}

\nt{\bf Proof}. Enough to show the first statement. Let $\ov{P}_U$
denote the complex frame bundle of $(X,-J)$. Since
$\si:(X,J)\lra(X,-J)$ is holomorphic, it induces a unique bundle
isomorphism ${\si}_*: P_U\lra\ov{P}_U$. Note that
$\ov{P}_U=P_U\times_{c} U(n)$, where ${c(u)=\ov{u}}$ is the
conjugation map on $U(n)\lra U(n)$. Thus one can take $\tau$ to be
${\si}_*$ composed with the conjugation $c: \ov{P}_U\lra {P}_U$.
 \qed

In particular, the lifting $\tau:P_{so}\to P_{so}$ is {\em not} the
induced map $\si_*:P_{so}\to P_{so}$, since the latter is a
$SO(2n)$-bundle isomorphism.

 Recall the embedding $\gamma: U(n)\lra Spin^c(2n)$ can be
defined as follows (cf. Lawson-Michelson \cite{lm}): if $u\in
U(n)$ is diagonalized as
 $$\text{diag}\{e^{{\bf i}\theta_1}, \cdots, e^{{\bf
i}\theta_n}\}$$ under a complex basis $(\epsilon_1,
\cdots,\epsilon_n)$ of ${\bf C}^n$, then $\gamma(u)\in
Spin(2n)\times_{\pm} U(1)=Spin^c(2n)$ is
$$\prod_k (\cos\frac{\theta_k}{2}+
\sin\frac{\theta_k}{2}\cdot\epsilon_k\cdot J\epsilon_k)\times
e^{\frac{\bf i}{2}\sum_k\theta_k}.$$
 Since $\ov{u}$ is diagonalized as $\{e^{-{\bf i}\theta_1},
 \cdots, e^{-{\bf i}\theta_n}\}$
under the complex basis $(\ov{\epsilon}_1,
\cdots,\ov{\epsilon}_n)$, $\gamma(\ov{u})$ is equal to
$$
\prod_k (\cos\frac{\theta_k}{2}-
\sin\frac{\theta_k}{2}\cdot\ov{\epsilon}_k\cdot
J\ov{\epsilon}_k)\times e^{-\frac{\bf i}{2}\sum_k\theta_k},$$
which is the same as
 $$\ov{\prod_k (\cos\frac{\theta_k}{2}+
\sin\frac{\theta_k}{2}\cdot\epsilon_k\cdot J\epsilon_k)}\times
e^{-\frac{\bf i}{2}\sum_k\theta_k},$$
 by noting that
 $\ov{J\epsilon_k}=\ov{\i\epsilon_k}=-J\ov{\epsilon}_k$.
 Here the conjugation on the $Spin(2n)$-factor is the restriction
 of that to the Clifford algebra $Cl({\bf R}^{2n})=Cl({\bf C}^n)$,
 namely the one generated by the standard conjugation on ${\bf
 C}$. This means that only by conjugating {\em both}
 factors $Spin(2n)$ and
 $U(1)$, we obtain the conjugation on $Spin^c(2n)$ which is
  compatible with the
 conjugation on $U(n)$ via the inclusion $\gamma$. (One
 should emphasize that the other conjugation on $Spin^c(2n)$ coming from
  the $U(1)$-factor alone, as usually considered, is
 not what is required here.) In general we call this kind of
 conjugation coming from both factors  a
 ``diagonal conjugation''.

Now $Cl({\bf R}^2)\otimes_{\bf R}{\bf C}$ is canonically
isomorphic to the matrix algebra ${\bf C}[2]$, and the diagonal
conjugation on $Cl({\bf R}^2)\otimes_{\bf R}{\bf C}$ -- from both
factors, is compatible with the usual entry wise conjugation on
${\bf C}[2]$ under this isomorphism. Using the periodicity
$Cl({\bf R}^{2n})\otimes{\bf C}\cong [Cl({\bf R}^2)]\otimes{\bf
C}]^{\otimes n}= {\bf C}[2^n]$, it is not hard to check that the
diagonal conjugation on $Cl({\bf R}^{2n})\otimes{\bf C}$ is
compatible with the standard entry wise complex conjugation on
${\bf C}[2^n]$.

Next consider the standard complex spin representation
$$Spin^c(2n)\hookrightarrow Cl({\bf R}^{2n})\otimes{\bf
C}={\bf C}[2^n]\hookrightarrow GL_c(V),$$
 where $V$ is a complex vector space of dimension $2^n$. Then the
 diagonal conjugation on $Spin^c(2n)$ is compatible with the
 complex conjugation on $V$.

Recall that  the canonical spin$^c$ bundle is
$P_{sp}=P_U\times_{\gamma}Spin^c(2n)$ and the associated spinor
bundle is $S=P_{sp}\times V$.

\begin{pro}\label{spinl}
There is a canonical lifting $\tau:P_{sp}\to P_{sp}$ of
$\si$, which satisfies
$$\tau(pg)=\tau(p)\ov{g},$$
where $p\in P_{sp}$ and $\ov{g}$ signifies the diagonal
conjugation on $Spin^c(2n)$.

The induced lifting on $S$ (still denoted by $\tau$) is fiberwise
complex anti-linear and compatible with the complex Clifford
multiplication $(T^*X\otimes{\bf C})\times S\to S$, where
$\si^*:T^*X\to T^*X$ should be extended as anti-linear map on the
complexification $T^*X\otimes{\bf C}$.
\end{pro}

\nt{\bf Proof}. The lifting on $P_{sp}$ is induced from the one
given in Lemma \ref{u}. The compatibility holds because the
lifting on $S$ is constructed via the spin$^c$-principal bundle
and the diagonal conjugation on $Spin^c(2n)$ is compatible with
the conjugation on $V$ as discussed above.\qed

Note that the determinant line bundle of $P_{sp}$ also carries a
natural anti-linear lifting of $\si$, since the isomorphism
 $U(1)\cong U(1)/\pm1$ preserves the complex conjugation. In fact,
 $\det P_{sp}=K^{-1}$ (anti-canonical bundle of $J$), which
 certainly has an anti-linear lifting. More generally, for any
 line bundle $L'\to X$ such that $\si^*c_1(L')=-c_1(L')$, the
 corresponding $Spin^c$ bundle $S_{L'}=S\otimes L'$ (with
 determinant $K^{-1}\otimes(L')^2$) has a canonical anti-linear
 lifting, compatible with the Clifford multiplication.

Consider now an arbitrary real vector space $W$  with an almost complex
structure $J$.
Given any linear map $\si:W\to W$ such that $\si\circ J=-J\circ\si$, we
 extend it s{\em anti-linearly}
on the complexification, $\tilde{\si}:W\otimes 
{\bf C}\to W\otimes {\bf C}$, so that
$\tilde{\si}(w\otimes c)=\si(w)\otimes\bar{c}$. This 
contrasts with the usual
{\em linear} extension of $J$ on the complexification, 
and is required by the following lemma.

\begin{lem}\label{ext} The extension $\tilde{\si}$ preserves
 the decomposition $W\otimes{\bf C}=W^{1,0}\oplus W^{0,1}$
of the $\pm{\bf i}$-eigen spaces of $J$.
\end{lem}

\nt{\bf Proof}. Take any $w\in W^{1,0}$. Then $J(w)={\bf i}w$.
Since $J(\si(w))=-\si(J(w))=-\si({\bf i}w)={\bf i}\si(w)$, we have
$\si(w)\in W^{1,0}$. This establishes $\si: W^{1,0}\to W^{1,0}$.
The second summand is similar. \qed\\

\nt{\bf Remark.} Because of this lemma, from now on we will always take
the anti-linear extension of an anti-holomorphic involution $\si$ on the
complexification. We will also use $\si$ for the extension without the
tilde sign.

\begin{cor}\label{liff}
Identify the spinor bundle $S$ canonically with the cotangent bundle
$\Lambda^{0,*}X=
\oplus_r\La^{0,r}X$ of $(0,*)$-forms  as usual. Then the lifting $\tau$ on $S$
is equivalent to the anti-linear lifting $\si^*$ on $\Lambda^{0,*}X$.
\end{cor}

\nt{\bf Proof}. From the early discussion,
$$S=P_{sp}\times V=(P_U\times_\gamma Spin^c(2n))\times V
=P_U\times_\eta V,$$
where the composition $\eta: U(n)\to GL_c(V)$ of $\gamma$ with the
spin representation is the standard unitary representation on
$V=\Lambda^*{\bf C}^{n}$. Thus $S=\La^*T_cX$, where $T_cX$ is the
tangent bundle with almost complex structure $J$ and the wedge
product is taken over ${\bf C}$ fiberwisely. Clearly the lifting
$\tau$ on $S$ is equivalent to the lifting $\si_*: T_cX\to T_cX$
which is fiberwise anti-linear, since $\si_*\circ J=-J\circ\si_*$.

The natural identification of $S=\La^*T_cX$ with $\La^{0,*}X$ is
through
$$\begin{array}{ccc}
T_cX&\longrightarrow& \La^{0,1}X\\
v&\mapsto&v^*=(\bullet,v)
\end{array}$$
where the metric is used. Since $\si$ is an isometric involution,
$\si_*(v)\mapsto (\bullet,\si_*(v))=(\si_*(\bullet), v)=\si^*(v^*)$.
Hence $\si_*$ is equivalent to the lifting $\si^*$ on $\La^{0,1}X$.
That $\si^*$ is anti-linear follows from $\si^*\circ J=-J\circ\si^*$ and
$J=-{\bf i}$ on $\La^{0,1}X$.\qed

\subsection{Seiberg-Witten equations with real structures}\label{srw}

Now we specialize to the case of our interest, that  $(X, J)$ is a
Hermitian 4-dimensional almost complex manifold with an
isometric anti-holomorphic
involution $\si$.

\nt{\bf Convention.} When no confusion is possible,
 we will often use $\ov{w}$ for $\si(w), \si_*(w), \si^*(w)$ or
more generally for $\tau(w)$,
where $\tau$ is any induced map by $\si$. This is  convenient and makes
sense since the maps are often anti-linear.

For example,  $\ov{{\bf i}\al}=-{\bf i}\ov{\al}$ interprets
conveniently the formula $\si^*({\bf i}\al)=-{\bf i}\si^*(\al)$
with ${\bf i}\al\in \Om^*\otimes {\bf C}$.
In particular, if $\i\al=F$ is the curvature 2-form of a unitary connection on
a complex line bundle over $X$, then
$\ov{F}=-\i\ov{\al}$.
The appearance of the $-$ sign here will save a lot of $-$ signs
elsewhere.

\begin{lem}\label{Llift}
 Suppose $L\lra X$ is a complex line bundle and $\tau: L\to L$
is an anti-linear lifting of $\si$ (so $\si^* c_1(L)=-c_1(L)$).
For any unitary connection $A$ on $L$ and its pull-back
$\ov{A}=\tau^*A$, their curvatures satisfy $F_{\ov{A}}=\ov{{F}_A}$.

\end{lem}

\nt{\bf Proof}.  One can prove the lemma by direct calculations on
the local connection matrices under a gauge. More convenient is to
use the corresponding principal bundle $P$ of $L$. Then the
lifting $\tau: P\lra P$ satisfies $\tau(pg)=\tau(p)\ov{g}$ for
$p\in P, g\in U(1)$. The connection 1-form $\omega_A$ is globally
defined on $P$ with values in the Lie algebra ${\bf i}{\bf R}$ of
$U(1)$. Since the conjugation $g\mapsto \ov{g}$ induces
 the map $\xi\mapsto -\xi$ on ${\bf
i}{\bf R}$, which is compatible with the anti-complex linear
extension of $\si$, the connection form
of $\ov{A}$ is
$\omega_{\ov{A}}=\ov{\omega_A}$. It follows that
$F_{\ov{A}}=\ov{{F}_A}$. \qed

 Thus if $A$ is equivariant under
$\tau$, then its curvature obeys $F_A=\ov{{F}_A}$.

Consider  the canonical spin$^c$ bundle $S$ of $(X,J)$, with
determinant bundle $L=K^{-1}$. By Proposition \ref{spinl}, we have
a canonical anti-linear lifting $\tau$ on $S$ and $L$.

Under the previous remark, for a section
$\Phi\in\Ga(S^+)$, $\ov{\Phi}$ is the pull-back section
$\tau^*\Phi:=\tau^{-1}\circ\Phi\circ\si$. Similarly,
  $\ov{h}={\si}^*(h)$ if $h\in\i\Om^2_+$. The induced action
of $\si$ on the gauge group $\mc{G}$ is
$ g\mapsto\ov{g}$ where $\ov{g}(x)=\ov{g(\ov{x})}$, i.e.
$\ov{g(\si(x))}$, the long over line
being the conjugation on $S^1$. With these actions understood, we
have the following:

\begin{pro}\label{rsw} (i) The gauge transformation
$\mc{G}\times\mc{C}\to\mc{C}$ is $\si$-equivariant, where
$\mc{C}=\mc{A}\times\Ga(S^+)$ is the configuration space. Hence
the quotient space $\mc{B}=\mc{C}/\mc{G}$ has an induced
involution $\si$.

(ii) The SW function
$\mc{F}:\mc{C}\times\i\Om^2_+\to\i\Om^2_+\times\Ga(S^-)$ is
equivariant also. Hence $(A,\Phi)$ is a SW solution with respect
to $h$ iff $(\ov{A},\ov{\Phi})$ is a SW solution with respect to
$\ov{h}$.

(iii) The projection $\pi:\mc{C}\times\i\Om^2_+\to\i\Om^2_+$ is
equivariant, so is $\pi:{\bf  M} \to\i\Om^2_+$ after dividing
gauge transformations.
 \end{pro}

\nt{\bf Proof}. The statements follow from Proposition \ref{spinl}
coupled with Lemma \ref{Llift}. \qed

 Note that we may also prove the proposition using Corollary \ref{liff}.
 \\

Proposition \ref{rsw} can be obviously extended from the canonical
spin$^c$ structure $P_{sp}$ to a general one:

\begin{pro}\label{rsw1} If a principal $U(1)$-bundle $\xi$ has
an anti-linear lifting of $\si$, then all three parts of \ref{rsw}
remains to be true for the twisted spin$^c$ structure $P_{sp}\times\xi$
  of $P_{sp}$ by $\xi$.
\end{pro}

\nt{\bf Remark.} It is important to point out that the
parameterized moduli space ${\bf M}$ is not a complex or almost
complex manifold, partly because $\i\Om^2_+$ is not so.
Furthermore, the fibers of $\pi$ do not have any obvious complex
structure, except the un-perturbed moduli space $\pi^{-1}(0)$ on a
K\"{a}hler surface $X$.
 Nonetheless, it is convenient to say $\pi:{\bf M}\to \i\Om^2_+$ is {\em
real} which simply is taken to mean that $\pi$ is $\si$-equivariant.
By the same token, even though $\i\Om^2_+$ is not a complex space,
we still call $(\i\Om^2_+)_\R:=\mbox{Fix}(\si:\i\Om^2_+\to \i\Om^2_+)$
the real space. Note that under the convention above, $(\i\Om^2_+)_\R$
consists of  $\i\cdot (\si$-anti-invariant smooth
forms), namely $\i\al$ where $\ov{\al}=-\al\in \Om^2_+$.
 Similar remark applies
to $(\i\Om^r)_\R$ of other degrees. This is consistent with
$\ov{c_1(L)}=\si^*c_1(L)=-c_1(L)$.

On various occasions  it will be useful to define real liftings in
a topological way, irrespective of any almost complex structure on
$X$. Proposition \ref{spinl} motivates the following:

\begin{df}\label{tsi}
 Let $X$ be a smooth manifold of dimension $2n$  and $\si:X\to X$ a
smooth involution that admits a conjugate lifting on the frame
bundle $P_{so}$, $\si(pv)=\si(p)\ov{v}$, where $\ov{v}=TvT^{-1}$ as
in (\ref{dt}).

(i) A spin structure on $X$ is called {\em real compatible with} 
$\si$ if
the $Spin(2n)$-bundle $P_s$ admits a conjugate lifting $\tau$ of
$\si$, namely $\tau(pg)=\tau(p)\ov{g}$, where for ${g}\in Spin(2n)$,
$\ov{g}$ is the restriction of the complex conjugation from
$Cl({\bf C}^n)$.

(ii) Similarly a spin$^c$ structure on $X$ is real compatible with $\si$
if its principal $Spin^c(an)$-bundle $P$ admits a conjugate
 lifting $\tau$, $\tau(pg)=\tau(p)\ov{g}$, where $\ov{g}$ is
the diagonal conjugation of $g\in Spin^c(2n)$.
\end{df}

In both cases we will also call $\tau$ (topological) 
{\em real liftings}.
In terms of the spinor bundle $S$, $\tau$ leads to an
anti-linear involution lifting on $S$ which is compatible
with the Clifford
 multiplication on $T^*X\otimes{\bf
 C}\hookrightarrow\text{End}_c(S)$. As before, $\si$ should be
 extended as an anti-linear map on $T^*X\otimes{\bf C}$ in order
 to have this compatibility. In particular for $\dim X=4$ and a
compatible spin$^c$ structure $P$,
 with the  same induced action on the gauge group $\mc{G}$ and Lemma
 \ref{Llift} as in the previous section, Proposition \ref{rsw}
 carries over to the new set-up. In particular, $\mc{B}$ inherits
 an involution, and $(A,\Phi)$ is a SW solution with perturbation $h$
 iff $(\ov{A},\ov{\Phi})$ is with perturbation $\ov{h}$.

\section{Configuration spaces and determinant bundles}\label{dlin}

In this and next sections, to be definitive,
 we  focus on a K\"{a}hler surface
$(X,\om, J)$ that has an
 isometric real structure $\si$, thus $\si^*\om=-\om, 
\si^*\circ J=-J\circ\si^*$. We will indicate when appropriate 
that
many results below either remain to be true
(for example those in Subsection \ref{rcm}) or can be modified
suitably for an  almost complex or symplectic  manifolds.

Suppose that $S=S^+\oplus S^- \to X$ is a spin$^c$ structure
admitting a real lifting of $\si$ (cf. Proposition \ref{rsw1}).
Let  $L=\det S^+$ be the determinant bundle of the spin$^c$ structure.

\subsection{The real configuration and moduli spaces}\label{rcm}

 Set $\mc{C}^*=\mc{A}(L)\times(\Ga(S^+)\backslash 0)$. By Proposition
 \ref{rsw1}, the induced involutions on $\mc{C}^*, \mc{G}$, namely
$(A,\Phi)\mapsto (\ov{A}, \ov{\Phi}), g\mapsto\ov{g}$, are compatible:
$$\ov{g\cdot(A,\Phi)}=\ov{g}(\ov{A},\ov{\Phi}).$$
Thus we have the further induced involution $\si$ on  the 
configuration space
 $\mc{B}^*=\mc{C}^*/\mc{G}$ and hence the
{\em fixed configuration space}
$$\mc{B}^{*\si}=\mbox{Fix}(\si: \mc{B}^*\to \mc{B}^*)\subset \mc{B}^*.$$
Moreover, we have the {\em real configuration space} defined as
$$\mc{B}^*_\R=\mc{C}^*_\R/\mc{G}_\R$$
namely the set of real points of $\mc{C}^*$ modulo the real gauge group.
It follows essentially from the compatibility above and the 
freeness of the
$\mc{G}$ action on $\mc{C}^*$ that the   natural map
$[(A,\Phi)]_\R\mapsto [(A,\Phi)]$ gives rise to a
inclusion $\mc{B}^*_\R\hookrightarrow\mc{B}^{*\si}$.
 (In the appendix, we  organize and state the results for the general
set-up.) In this paper, we will be mainly interested in the space
$\mc{B}^*_\R$
and its subspace of real Seiberg-Witten solutions.

From the standard Seiberg-Witten theory, e.g. the book \cite{m},
the gauge group
$\mc{G}$ is naturally homotopic to $S^1\times H^1$, where
$H^1=H^1(X,{\bf Z})$.
As the classifying
space of the group $\mc{G}$,  $\mc{B}^*$ is weakly homotopic to
${\bf C}{\bf P}^{\infty}\times K(H^1, 1)$. Since $\si$ induces 
the standard
conjugation on the $S^1$-factor through $\mc{G}$, the induced action
on ${\bf C}{\bf P}^{\infty}$ is also the conjugation. Hence taking
fixed points
on both sides, we have
$$\mc{B}^{*\si}\sim{\bf R}{\bf P}^{\infty}\times K(H^1, 1)^\si.$$
An similar argument will give the following result for the 
weak homotopy type of $\mc{B}^*_\R$.

\begin{pro}\label{hom} There is a natural weak homotopy equivalence:
$$\mc{B}^*_\R\sim {\bf R}{\bf P}^{\infty}\times K(H^1_\R, 1),$$
where $H^1_\R=H^1(X,{\bf Z})^\si$.
\end{pro}

\nt{\bf Proof}. The real constant gauges form a subgroup:
 ${\bf Z}_2\subset\mc{G}_\R$. For the quotient group, there is a natural
bijection
$\pi_0(\mc{G}_\R/{\bf Z}_2)\to H^1_\R$ given by 
$\xi\mapsto \rho_\xi$, where
$\rho_\xi\in H^1(X,{\bf Z})^\si\subset H^1(X, {\bf R})^\si$ 
is defined as
 $\rho_\xi=[g^{-1}dg]$ for a gauge $g\in\xi$ such that $g^{-1}dg$ is a
real harmonic 1-form.
It follows that $\mc{G}_\R$ is homotopic to 
${\bf Z}_2\times H^1_\R$ and
the classifying space $B\mc{G}_\R$ of $\mc{G}_\R$ is
weakly homotopic
to ${\bf R}{\bf P}^{\infty}\times K(H^1_\R, 1)$. Since the real part
$\mc{C}^*_\R$ is clearly contractible just as $\mc{C}^*$ is, the real
configuration space $\mc{B}^*_\R=\mc{C}^*_\R/\mc{G}_\R$ is 
weakly homotopic to
$B\mc{G}_\R$ hence to ${\bf R}{\bf P}^{\infty}\times K(H^1_\R, 1)$.
\qed\\

\nt{\bf Remarks.} (1) The generator in $H^2(\mc{B}^*,{\bf Z})$ 
that comes from
the ${\bf C}{\bf P}^{\infty}$ factor restricts to a 2-torsion in
$H^2(\mc{B}^*_\R,{\bf Z})$. In fact the resulting complex line bundle
on $\mc{B}^*_\R$ is the complexification of the real line bundle
corresponding to the generator in $ H^1(\mc{B}^*_\R,{\bf Z}_2)$
that comes from the ${\bf R}{\bf P}^{\infty}$ factor.

(2) By Proposition \ref{hom}, $\mc{B}^*_\R$ is connected; in contrast,
the fixed configuration space $\mc{B}^{*\si}$ is disconnected and
contains $\mc{B}^*_\R$ as a connected component.

At a real point $(A,\Phi)\in\mc{C}^*_\R$, the tangent space is
$$T_{A,\Phi}\mc{C}^*_\R=(\i\Om^1)_{\R}\oplus\Ga(S^+)_\R,$$
where the subscript $\mbox{R}$ indicates the invariant subspaces
under the extended $\si$-action.
Linearizing the $\mc{G}_\R$ action on $\mc{C}^*_\R$, we have
$$D^1_\R:(\i \Omega^0)_\R\to (\i\Omega^1)_\R\oplus\Ga(S^+)_\R$$
as the restriction of $D^1$ from the complex (\ref{complex}). Thus
the tangent space $T_{[A,\Phi]}\mc{B}^*_\R=\ker(D^1_\R)^*$,
from which one sees that $\mc{B}^*_\R$ is an open subspace 
of $\mc{B}^{*\si}$.

The relevant complex for the real parameterized moduli space is
\begin{equation}\label{rcomplex}
 \begin{array}{ccccc}
&D^1_\R&&D^2_\R&\\
0\lra(\i \Omega^0)_\R&\lra&(\i\Omega^1)_\R\oplus\Ga(S^+)_\R&\lra&
(\i\Omega^2_+)_\R\oplus\Ga(S^-)_\R\lra0,
\end{array}
\end{equation}
by restricting the complex (\ref{complex}) to the real spaces.

By Proposition \ref{rsw1}, the Seiberg-Witten functional restricts to
the real spaces:
\begin{equation}\label{swfr}
\begin{array}{cccc}
\mc{F}_\R:&\mc{A}_\R\times\Ga(S^+)_\R\times(\i\Om^2_+)_\R&\lra&
(\i\Om^2_+)_\R\times\Ga(S^-)_\R.
\end{array}\end{equation}
At a real point $(A,\Phi,h)$, the differential
$$D\mc{F}_\R:\;
(\i\Om^1)_\R\oplus\Ga(S^+)_\R\oplus(\i\Om^2_+)_\R \lra(\i\Om^2_+)_\R
\oplus\Ga(S^-)_\R$$
 is  $D\mc{F}_\R(a,\phi,k)=D^2_\R(a,\phi)-(k,0)$.
  The usual proof of the transversality theorem
can be adapted easily to show that $0$ is
a regular value of $\mc{F}_\R$ when restricted to irreducibles, hence
$\mc{F}^{-1}_\R(0)\cap \{\Phi\not=0\}$ is a smooth
Banach manifold. The tangent space at the real point $(A,\Phi,h)$ is
$T\mc{F}^{-1}_\R(0)^*=\ker D\mc{F}_\R$.

Dividing by real gauge transformations, we have the {\em
parameterized real moduli spaces}: ${\bf  M}_\R
=\mc{F}^{-1}_\R(0)/\mc{G}_\R$. The real version of Proposition
\ref{dpi} becomes:

\begin{pro}\label{dpir}
 (i) The tangent space of the parameterized irreducible real moduli
 space ${\bf  M}^*_\R $
 at a real point $[(A,\Phi,h)]$ is the following subspace of
$(\i\Om^1)_\R\oplus\Ga(S^+)_\R\oplus(\i\Om^2_+)_\R$:
 $$\mr{T}{\bf  M}^*_\R =\{(a,\phi,k)\mid D^2_\R(a,\phi)=(k,0) 
\text{ and }
(D^1_\R)^*(a,\phi)=0\}.$$

(ii) The differential $D\pi_\R$ of the projection map
$\pi_\R:{\bf  M} ^*_\R\to(\i\Om^2_+)_\R$ has the kernel and cokernel
canonically identified with:
$$\ker(D\pi_\R)={\bf H}^1_\R, \; \mr{coker}D\pi_\R={\bf H}^2_\R,$$
where ${\bf H}^1_\R, {\bf H}^2_\R$
are the cohomology of the complex (\ref{rcomplex}).
\end{pro}

 A single elliptic operator converted from the basic complex
(\ref{rcomplex}) is
\begin{equation}\label{rde}
 \de_\R=D^2_\R\oplus
(D^1_\R)^*:(\i\Om^1)_\R\oplus\Ga(S^+)_\R\lra [(\i\Om^2_+)_\R\oplus
\Ga(S^-)_\R]\oplus(\i\Om^0)_\R.
 \end{equation}
The real version of Corollary \ref{dee} is

\begin{cor}\label{rdee} There are natural isomorphisms
 \begin{equation}\label{rde0}
  \ker D\pi_\R\cong
 \ker\de_\R,\; \mr{coker}D\pi_\R\cong\mr{coker}\de_\R.
\end{equation}
\end{cor}

Thus the {\em orientation bundle}
$\det\pi_R=\bigwedge^{\mr{max}}\ker D\pi_\R\otimes
 (\bigwedge^{\mr{max}}\mr{coker}D\pi_\R)^*$ of the map $\pi_\R$
is naturally identified with the determinant bundle of $\de_{\mr{R}}$:
\begin{equation}\label{oo}
\det\pi_R=\bigwedge^{\mr{max}}\ker\de_{\mr{R}}\otimes
 (\bigwedge^{\mr{max}}\mr{coker}\de_{\mr{R}})^*
\end{equation}
 on ${\bf  M} ^*_\R$.
This is the reason why the latter bundle will play a prominent
role in the paper.

\subsection{The real determinant line bundle}\label{rdl}

At a point $(A,\Phi)\in\mc{C}^*$, let us decompose the operator
 $\de=\de_{A,\Phi}:\i\Om^1\oplus\Ga(S^+)
 \to[\i\Om^0\oplus\i\Om^2_+]\oplus\Ga(S^-)$
defined in (\ref{de}) as
 $\de=(\de^X\oplus\Dr_A)+\eta$ where
\begin{equation}\label{ff}
\de^X=(2d^*,d^+):\i\Om^1\to \i\Om^0\oplus\i\Om^2_+
 \end{equation}
  depends on $X$ only and
 $\eta=\eta_{_{\Phi}}$ is a zero-th order operator depending on
$\Phi$ only:
 $$\eta_{_{\Phi}}(a,\phi)=(-<\phi,\Phi>-Dq_{\Phi}(\phi))+\frac{1}{2}a
\cdot\Phi.$$
 Note that $\eta_{_{t\Phi}}=t\eta_{_{\Phi}}$; in particular
 $\eta_0$ is
the zero operator. If we set further $\de^L=\de^L_A=\de^X\oplus\Dr_A$,
then
$\de=\de^L+\eta$ so that $A, \Phi$ are separated in the two operators.

 \begin{pro}\label{dex}  There are suitable
 complex structures in the infinitely dimensional spaces
$\i\Om^1, \i\Om^0\oplus\i\Om^2_+$
 such that the extended
$\si$ actions are anti-holomorphic on these spaces and $\de^X$ is
  complex linear. Moreover, the numerical index
$\mr{ind}\de_{\mr{R}}$  is half of $\mr{ind}\de$
namely
$$ \mr{ind}\de_{\mr{R}}=\frac{1}{8}(c_1(L)^2-2e_X-3s_X).$$

 \end{pro}

 \nt{\bf Proof.}
The space $\i\Om^1\cong\i\Om^{0,1}$  has the induced complex
structure by $J$ under which $\si$ is anti-holomorphic.
Whereas there is a natural real linear isomorphism
 $\Om^0\oplus\Om^2_+\cong\Om^0\oplus\Om^0\cdot\om\oplus\Om^{0,2}$,
 the latter being isomorphic to $\Om^0_c\oplus\Om^{0,2}$
 by viewing $\om=\i$ on the complexification $\Om^0_c$.
Hence $\i\Om^0\oplus\i\Om^2_+$ inherits a complex
 structure, under which  $\si$ is anti-holomorphic in view of
 $\si^*\om=-\om$.
 Furthermore $\de^X$ is complex linear, since it is
 equivalent to $\ov{\partial}^*\oplus\ov{\partial}:\Om^{0,1}\to
 \Om^0_c\oplus\Om^{0,2}$ under the previous transformations.
(In the symplectic case, they are equivalent up to a 
zeroth order operator.)

At a real point $(A,\Phi)\in\mc{C}_\R$,
$\de^L=\de^X\oplus\Dr_A$ is complex linear and real with respect to
$\si$. Hence $\ker\de^L, \cok\de^L$ are complex vector spaces with
real structure,
and the real parts have half the dimensions, giving
$\mr{ind}\de^L_{\R}=\mr{ind}\de^L/2$. Since $\eta, \eta_{\R}$ are
zeroth order operators, the indices remain the same for
$\de=\de^L+\eta$ and $\de_{\R}=\de^L_{\R}+\eta_{\R}$. Thus
$\mr{ind}\de_{\R}=\mr{ind}\de/2$ holds.
\qed

With respect to the extended real structure, the previously 
defined fixed
point set $(\i\Om^1)_\R$ is now the true real part of $\i\Om^1$. 
Clearly
the real part of $\i\Om^0\oplus\i\Om^2_+$ is
$$[\i\Om^0\oplus\i\Om^2_+]_\R=(\i\Om^0)_\R\oplus(\i\Om^2_+)_\R,$$
where the summands are fixed point sets of $\si$ (which are 
not real parts).
Note that for when $\Phi\not=0$, $\eta$ is not a complex linear operator,
because of the
quadratic term $Dq_{\Phi}(\phi)$.  (But it is $\si$-equivariant
and so $\eta_{\R}$ is defined, as we have used above.) Thus unlike
$\ker\de^L, \cok\de^L$, the spaces
$\ker\de, \cok\de$ are not necessarily complex vector spaces.

By Proposition \ref{dex},  $\de_\R=\de_{A,\Phi;\R}$ is certainly
a Fredholm operator, which is
 parameterized by $(A,\Phi)\in\mc{C}_{\R}$.
As usual such a Fredholm family gives rise to
the (real) determinant line bundle
$$\det\mr{ind}\de_{\mr{R}}=\bigwedge^{\mr{max}}\ker\de_{\mr{R}}\otimes
 (\bigwedge^{\mr{max}}\mr{coker}\de_{\mr{R}})^*,$$
which descends to the real configuration space ${\mc B}^*_{\mr{R}}$,
since the real gauge group $\mc{G}_\R$ action lifts to the bundle.
We still denote the descended bundle by $\det\mr{ind}\de_{\mr{R}}\to
{\mc B}^*_{\mr{R}}$.
(In the almost complex surface case, $\de_\R$ is still Fredholm,
because $\ker\de_\R\subset \ker\de$ and $\cok\de_\R\subset\cok\de$
both are finite dimensional.
The second inclusion uses $\cok\de_\R=\ker\de^*_{\R}$, $\si$
is isometric, etc.) The bundle is an extension of $\det\pi_{\R}$ in view
of \ref{oo}.

 Since $\pi_{\R}$ is clearly a
Fredholm map, by Sard-Smale theorem, the regular values of
$\pi_{\R}$ form a dense subset of $(\i\Om^2_+)_{\R}$. For 
each regular value
$h$, the corresponding real moduli space $M_{\R}(h)=\pi^{-1}_{\R}(h)$
is a smooth manifold.
As in the usual situation, its
orientation bundle is the restriction of $\det\mr{ind}\de_{\R}$
 to $M_{\R}(h)\subset\mc{B}^*_{\R}$. However, the  orientation
  of $\det\mr{ind}\de_{\R}\to\mc{B}^*_{\R}$ is much
 more complicated in the current real case. Indeed we will see that
the bundle is in general non-orientable (i.e. non-trivial).

Let $H^i_\R(X,\R)$ denote the {\em real} De Rham cohomology group
with respect to $\si$, namely the space
of $\si$-invariant closed forms modulo $\si$-invariant exact forms.
While we define the {\em fixed}
cohomology $H^i(X,\R)^\si=\mbox{Fix}(\si^*:H^i(X,\R)\to H^i(X,\R))$.
Similarly introduce $H^+_\R(X,\R)$ and $H^+(X,\R)^\si$.
The following is a simple consequence of the classical Hodge theorem,
using only that $\si$ is isometric.

\begin{lem}\label{hog}
There are  natural isomorphisms
$$H^i_\R(X,\R)\cong H^i(X,\R)^\si,\;\;
H^+_\R(X,\R)\cong H^+(X,\R)^\si.$$
\end{lem}

\nt{\bf Proof}. To show $H^i_\R(X,\R)\cong H^i(X,\R)^\si$,
it is enough to show that the natural inclusion
$H^i_\R(X,\R)\hookrightarrow H^i(X,\R)^\si$ is surjective. Take any
fixed class in $H^i(X,\R)^\si$ and represent it by the harmonic
$i$-form $\al$. Hence $[\al]=[\ov{\al}]\in H^i(X,\R)^\si$. Since
$\si$ preserves the metric on $X$, $\ov{\al}$ is also harmonic.
As each class has a unique harmonic representative, one must have
$\al=\ov{\al}$; hence $[\al]$ comes from a class in $H^i_\R(X,\R)$
that is represented also by $\al$.

The second isomorphism can be proved similarly. \qed

To orient the determinant
$\det\mr{ind}\de_{\R}=\det\mr{ind}(\de^L_{\R}\oplus\eta_{_{\R}})
\to\mc{C}_\R$, as in the usual case, there are two slightly different
but equivalent approaches available. One is to deform 
the {\em fiber} of
  $\det\mr{ind}\de_{\R}$ at a given point $(A,\Phi)\in\mc{C}_\R$ by
deforming
the operator $\de_\R=\de_{A,\Phi;\R}$ in a family:
$$\de_{A,\Phi;\R}(t)=\de^L_{A;\R}+t\eta_{_{\Phi;\R}}, \; 
0\leq t\leq 1,$$
so obtaining the deformed fibers $\det\mr{ind}\de_{A,\Phi;\R}(t)$
over the same point $(A,\Phi)$. The other is to first deform the
{\em point} $(A,\Phi)$ itself in a path $(A,t\Phi), 0\leq t\leq
1$, and then restrict the bundle $\det\mr{ind}\de_{\R}$ to the
path, which becomes the bundle $\det\mr{ind}\de_{A,t\Phi; \R}\to
[0,1]$. The two approaches are interchangeable through the
relation
$$\de_{A,t\Phi;\R}=\de^L_{A;\R}+\eta_{_{t\Phi;\R}}=\de^L_{A;\R}+
t\eta_{_{\Phi;\R}}
=\de_{A,\Phi;\R}(t).$$
In the end, both approaches relate  $\xi$  with
the  bundle $\det\mr{ind}\de^L_{\R}\to\mc{C}_\R$ by taking $t=0$,
where we suppress again the subscript $A$ in the operator family.
We shall adapt the second approach in the argument below, which is
conceptually more clear.

As before, $\mc{A}_\R$ denotes the space of real connections on $L$.
Let $\mc{B}^L_\R=\mc{A}_\R/\mc{G}_\R$. Note the $\mc{G}_\R$ 
action has
a stabilizer $\pm1$ at every point in $\mc{A}_\R$, hence 
$\mc{G}_\R/\pm1$
acts freely and $\mc{B}^L_\R=\mc{A}_\R/(\mc{G}_\R/\pm1)$ is a smooth
Banach manifold. Clearly the natural forgetting map
\begin{equation}\label{proj}
p:\mc{B}^*_\R\to \mc{B}^L_\R, \;\; [A,\Phi]\mapsto [A]
\end{equation}
is a smooth map. The determinant bundle of the real Dirac operators
$$\Dr_{A,\R}:\Ga(S^+)_\R\to\Ga(S^-)_\R$$
parameterized by $A\in\mc{A}_\R$ obviously descends to
$\mc{B}^L_\R$, which we denote  by $\det\mr{ind}\Dr_{A,\R}$ or
simply $\det\mr{ind}\Dr_{\R}$.

\begin{thm}\label{q} Fix orientations on  $H^1(X,\R)^\si, 
H^+(X,\R)^\si$.
The determinant $\det\mr{ind}\de_{\R}\to\mc{B}^*_{\R}$
is isomorphic to the pull-back bundle, 
$p^*\det\mr{ind}\Dr_{A,\R}$, via
an isomorphism that is unique up to a positive continuous function. 
In other
words, the bundle
$$\det\mr{ind}\de_{\R}\otimes(p^*\det\mr{ind}\Dr_{A,\R})^{-1}\to
\mc{B}^*_{\R}$$
is orientable with a canonical orientation.
\end{thm}

\nt{\bf Proof}. Consider  the real full configuration space 
$\mc{B}_\R=\mc{C}_\R/\mc{G}_\R$, which is Hausdorff and contains
 $\mc{B}^L_\R$ as a singular submanifold.
Clearly one can extend the map $p$ in (\ref{proj}) to $\mc{B}_\R$
as a continuous map, which in turn  establishes a homotopy 
equivalence
$\mc{B}_\R\simeq \mc{B}^L_\R$ through the standard
deformation retraction
$$\Theta:\mc{B}_\R\times[0,1]\to \mc{B}_\R, \;\; ([A,\Phi], t)\mapsto
 [A,t\Phi].$$
Modulo the lifted $\mc{G}_\R$ action, the determinant bundle
$\xi\to \mc{C}_\R$ descends to a continuous line bundle
$\det\mr{ind}\de_\R\to\mc{B}_\R$. Using the retraction $\Theta$,
one sees that $\det\mr{ind}\de_\R$ is isomorphic to the pull-back
of $\det\mr{ind}\de_\R|_{\mc{B}^L_\R}=\det\mr{ind}\de^L_\R$, where
$\de^L_\R=\de^X_\R\oplus\Dr_{A;\R}$ is understood to be
parameterized by $[A]\in \mc{B}^L_\R$. Indeed the isomorphism can
be obtained by deforming the points in $\mc{B}_\R$ as follows:
take the fibers $f, f'$ of $\det\mr{ind}\de_\R$ over an arbitrary
point $[A,\Phi]$ and its projection $[A,0]$. Of course the
topological line bundle $\det\mr{ind}\de_\R$ is trivial along the
path $[A,t\Phi], 0\leq t\leq 1$. Any trivialization gives rise to
an isomorphism between $f$ and $f'$. Moreover the isomorphisms
obtained through different trivializations differ by positive
constants. (So the correspondence of the orientations of $f,f'$ is
independent of the trivializations. Alternatively this orientation
correspondence can be obtained through the two connected
components of the set $ \det\mr{ind}\de_\R \backslash \{
0-\mbox{section}\}$ over the path.) Since the path depends
continuously on the point $[A,\Phi]$, one can choose a global
bundle isomorphism $\det\mr{ind}\de_\R\cong
p^*\det\mr{ind}\de^L_\R$ on $\mc{B}_\R$, which is unique up to a
positive continuous function. In particular the bundle
\begin{equation}\label{yyt}
\det\mr{ind}\de_{\R}\otimes(p^*\det\mr{ind}\de^L_\R)^{-1}\to
\mc{B}^*_{\R}
\end{equation}
is orientable with a canonical orientation.

Next we examine the bundle $\det\mr{ind}\de^L_\R\to\mc{B}^L_\R$.
Recall
 $$\det\mr{ind}\de_{\R}^L=
 \det\mr{ind}\de_{\R}^X\otimes\det\mr{ind}\Dr_{A,\R},$$
and $\det\mr{ind}\de_{\R}^X=\det\mr{ind}(2d^*_\R\oplus d^+_\R)$ is
a constant 1-dimensional vector space independent of $A$. An
orientation of $\det\mr{ind}\de_{\R}^X$ is determined by
orientations of the cohomology groups of the complex
\begin{equation}\label{rrcomplex}
 \begin{array}{ccccc}
&2d_\R&&d^+_\R&\\
0\lra(\i \Omega^0)_\R&\lra&(\i\Omega^1)_\R&\lra&
(\i\Omega^2_+)_\R\lra0.
\end{array}
\end{equation}

Since $\ker\de^X, \cok\de^X$  are  complex vector spaces  with natural
orientations by Proposition \ref{dex}, the orientations of 
the cohomology
groups of (\ref{rrcomplex}) are determined by those of the
``imaginary part'' complex
\begin{equation}\label{rrmplex}
 \begin{array}{ccccc}
&2d_\R&&d^+_\R&\\
0\lra(\Omega^0)_\R&\lra&(\Omega^1)_\R&\lra&
(\Omega^2_+)_\R\lra0.
\end{array}
\end{equation}
By Lemma \ref{hog}, the cohomology of the last complex are
isomorphic to $H^0(X,\R)^\si, H^1(X,\R)^\si$, 
$H^+(X,\R)^\si$.
Hence any orientations on  $H^1(X,\R)^\si$, $H^+(X,\R)^\si$
determine the isomorphism $\det\mr{ind}\de_{\R}^L\cong
 \det\mr{ind}\de_{\R}^X$. The theorem is proved by coupling with
(\ref{yyt}) above.
\qed

Note that $\mc{B}_\R^*$ is not homotopic to $\mc{B}^L_\R$, as
$\mc{B}_\R^*$ and $\mc{B}_\R$ are not homotopic. The latter is so
in spite that the complement $\mc{B}^L_\R$ has infinite codimensions.
 For example, the generator of $H^1(\mc{B}_\R^*,{\bf Z}_2)$
 that comes from the
${\bf R}{\bf P}^\infty$ factor according to Proposition \ref{hom}
does not extend over $ \mc{B}^L_\R$, since it restricts non-trivially on
the link of $\mc{B}^L_\R$ in $\mc{B}_\R$.

Unlike the complex Dirac operator family $\Dr_A$, the real family
$\Dr_{A,\R}$ in general produces non-orientable determinant bundle
$\det\mr{ind}\Dr_{A,\R}$.

\subsection{The real universal bundle}\label{ru}

We consider here more carefully the various universal bundles that
are related to
our index bundles in the previous subsections. First recall a
 universal complex line bundle
$\mc{L}\to\mc{B}^*\times X$ can be defined as the quotient bundle of
$\pi^*L\to \mc{C}^*\times X$ under the lifted
$\mc{G}$ action, where $\pi: \mc{C}^*\times X\to X$ is the projection
on the second factor.   One may also  define a universal bundle
on $\mc{B}^L\times X$ but the construction needs to be modified:
the $\mc{G}$ action on the space $\mc{A}$ of connections is not free,
so the quotient bundle of $\pi^*L\to\mc{A}\times X$ is undefined.
To overcome the problem, one needs to use the based connection space
and pull back the universal bundle constructed there. More precisely
choose any base point $x_0\in X$ and set $\mc{B}^L_0=\mc{A}/\mc{G}(x_0)$
where the based gauge group $\mc{G}(x_0)=\{g\in \mc{G}| g(x_0)=1\}$
acts freely. Thus the above construction yields again a universal
bundle $\mathbb{L}\to \mc{B}^L_0\times X$. On the other hand, there is
a natural identification $\mc{B}^L=\mc{B}^L_0$, through which one has the
universal bundle $\mathbb{L}\to \mc{B}^L\times X$ as a carry-over.

It is interesting to observe that the pull-back bundle
$\wt{p}^*\mathbb{L}$ is not isomorphic to $\mc{L}$, where
$\wt{p}:\mc{B}^*\times X\to \mc{B}^L\times X$ is the forgetting map:
$([A,\Phi], x)\mapsto ([A],x)$. Such a discrepancy originates from the
above varied construction of $\mb{L}$. In fact, by construction
 $\mathbb{L}$ restricts to a trivial bundle $\mathbb{L}_{x_0}$
 on the slice $\mc{B}^L_0\times\{x_0\}$, hence the pull-back
$\wt{p}^*\mathbb{L}_{x_0}\to\mc{B}^*\times\{x_0\}$ is trivial as well.
However the restriction $\mc{L}_{x_0}\to
\mc{B}^*\times\{x_0\}$ of $\mc{L}$ is non-trivial, since
$\mc{L}_{x_0}$ is the quotient of the trivial bundle
$\mathbb{L}_{x_0}'\to \mc{B}^*_0\times\{x_0\}$ under a free $S^1$
action. Here $\mc{B}^*_0=\mc{C}^*/\mc{G}(x_0)$ is the based
irreducible configuration space and $\mathbb{L}'\to \mc{B}^*_0
\times X$ is the universal bundle, constructed similarly as $\mb{L}$.
 There is a natural $S^1$ action
on $\mc{B}^*_0$ with quotient $\mc{B}^*$ and the principal circle
bundle $\mc{B}^*_0\to\mc{B}^*$ (the based point fibration)
associates exactly the vector bundle
$\mc{L}_{x_0}$. In other words, $c_1(\mc{L}_{x_0})$ is the
generator of $H^2(\mc{B}^*)$ from the ${\bf C}{\bf P}^\infty$-factor.

The real structure on $\mc{B}^*\times X$ lifts to an anti-linear
isomorphism on $\mc{L}$. It follows that one has a real
line bundle on $\mc{B}^{*\si}\times X_\R$ by restricting to
the fixed points.
Since $\mc{B}^*_\R\subset \mc{B}^{*\si}$, a further restriction gives us
the anticipated real universal line bundle
$\mc{L}_\R\to \mc{B}^{*}_\R\times X_\R$.

Using the Stiefel-Whitney class $w_1(\mc{L}_{\R})\in
H^1(\mc{B}^*_{\R}\times X_{\R},{\bf Z}_2)$ and the slant product
one defines a map,
$$\nu=w_1(\mc{L}_{\R})/: H_0(X_{\R},{\bf Z}_2)\to
H^1(\mc{B}^*_{\R},{\bf Z}_2).$$ This is in addition to the usual
map $\mu: H_0(X,{\bf Z})\to H^2(\mc{B}^*,{\bf Z})$ using the slant
product with $c_1(\mc{L})$. To make things less mysterious, let
$\mc{L}_{\mr{R},x_0}$ denote the restriction of $\mc{L}_{\R}$ to
 $\mc{B}^*_{\R}\times\{x_0\}$ where $x_0\in X_{\R}$. Then
$\nu(x_0)=w_1(\mc{L}_{\R,x_0})$, much like
$\mu(x_0)=c_1(\mc{L}_{x_0})$. Clearly the restriction of $\mu(x_0)$
 to $\mc{B}^*_\R$ is the complexification of $\nu(x)$, if they both
are viewed as bundles. In the end the classes $\mu(x_0), \nu(x_0)$
both are independent
of the point $x_0\in X_\R$, which we will simply call $\mu, \nu$, 
since
they respectively come from the ${\bf C}{\bf P}^\infty, 
{\bf R}{\bf P}^\infty
$ factors of $\mc{B}^*, \mc{B}^*_{\R}$.

By analogous constructions, one has the universal spin$^c$ bundle
$$\mb{S}=\mb{S}^{+}\oplus\mb{S}^{-}\to\mc{B}^L\times X $$
with  $\det\mb{S}^{+}=\mb{L}$. The last  bundle
 carries a tautological connection in the $X$ direction.
As a consequence, one obtains the virtual index bundle
$ \mr{ind}\Dr_{A}\in K(\mc{B}^L)$ and its real version that was used
in \ref{rdl}. The standard Atiyah-Singer family index theorem
can be applied to calculate the Chern character
ch$(\mr{ind}\Dr_{A})\in H^*(\mc{B}^L)$.

Return to the map $p:\mc{B}^*_\R\to \mc{B}^L_\R$ in (\ref{proj}).
This is a smooth fibration with
fibers $(\Ga(W^+)_\R-\{0\})/\pm1$ homotopic to ${\bf R}{\bf P}^\infty$.
It will be useful to settle the question whether the determinant bundle
$\det\mr{ind}\de_{\R}\to \mc{B}^*_\R$ can be isomorphic
to the bundle $\mc{L}_{\R,x_0}\to \mc{B}^*_\R$:

\begin{pro} The bundles $\mc{L}_{\R,x_0}$ and $\det\mr{ind}\de_{\R}$
 are never isomorphic. In other
words, $\nu\not=w_1( \det\mr{ind}\de_{\R})$.
 \end{pro}

\nt{\bf Proof.} As we have seen,
on each fiber of $p$, the class $\nu=w_1(\mc{L}_{\R,x_0})$ 
restricts to
the generator of $H^1({\bf R}{\bf P}^\infty,{\bf Z}_2)$. On the other
hand, by Theorem \ref{q},
$\det\mr{ind}\de_{\R}\cong p^*\det\mr{ind}\Dr_{A,\R}$ is a pull-back
bundle.
Hence $\det\mr{ind}\de_{\R}$ restricts trivially on fibers of $p$
and can not be isomorphic to $\mc{L}_{\R,x_0}$ as a result.
 \qed

\section{Real Seiberg-Witten invariants in orientable cases}\label{rss}

In \ref{rcm} and \ref{rdl}, we introduced the projection
 $\pi_\R: {\bf M}^*_\R\to (\i\Om^2_+)_\R$ from the parameterized
irreducible real moduli space. This is a Fredholm map. So by the
Sard-Smale theorem, for a generic perturbation $h\in(\i\Om^2_+)_\R$,
the real moduli space $M^*_\R(h)=\pi^{-1}_\R(h)$ is a smooth manifold
of dimension
$$m=\frac{1}{8}(c_1(L)^2-2e_X-3s_X).$$
(See Proposition \ref{dex}.) As in the standard case,
the same kind of a priori estimates can be
applied to real solution pairs $(A,\Phi)\in \mc{B}_\R$ (simply by
restriction) to show that each real moduli $M^*_\R(h)$ is compact,
provided that $h$ stays away from the real reducible wall
$$W_\R=\i c^+ + \mr{Im}d^+_{\R}\subset (\i\Om^2_+)_\R.$$
Here $c^+$ is the unique ($\si$ anti-invariant)
 self dual harmonic 2-form representing $c_1(L)$
and $d^+_{\R}:(\i\Om^1)_\R\to (\i\Om^2_+)_\R$ as before.
Note that $W_\R$ is an affine subspace  of codimension
$$b^+_{\R}:=\dim H^+_\R(X,\i\R)=\dim H^+(X,\R)^{-},$$
where the superscript $-$ indicates the $\si$ anti-invariant part
is used. For our K\"{a}hler manifold $X$ case, one can apply the
Hodge decomposition to show that $b^+_{\R}=1+p_g$, with $p_g$ the
geometric genus of $X$. Hence $b^+_{\R}>1$ iff $b^+=1+2p_g>1$.

Thus we have at least a ${\bf Z}_2$-fundamental class $[M^*_\R(h)]
\in H_m(M,{\bf Z}_2)$ for each generic perturbation $h\notin W_\R$.
Hence we can make the following definition.

\begin{df}\label{z2sw}
Suppose $\si$ is a real structure on a K\"{a}hler manifold $X$ and
$\xi=S^+\oplus S^-$ is a spin$^c$ structure on $X$, admitting a
real lifting of $\si$. One defines the ${\bf Z}_2$-valued real
Seiberg-Witten invariant to be the paring
$$SW_\R(\xi)=<[M^*_\R(h)], \nu\cup\cdots\cup\nu>,$$
where the cup product is taken $m$-times and $\nu\in H^1(\mc{B}^*,
{\bf Z}_2)$ as before. If $b^+_\R>1$ i.e. $b^+>1$, then $SW_\R(\xi)$
is independent of $h$. Otherwise it is well-defined in each of the two
chambers of $(\i\Om^2_+)_\R-W_\R$.
\end{df}

However, in view of the following result,
 such real Seiberg-Witten invariants are of limited usage
in most situations.

\begin{pro}
(i) When $X$ is of general type and $b^+>1$, the invariant
$SW_\R(\xi)$ is trivial unless $m=0$.

(ii) If $m=0$ (but for any $X$),  $SW_\R(\xi)$ is the mod 2
reduction of the
ordinary Seiberg-Witten invariant $SW(\xi)$.
\end{pro}

\nt{\bf Proof.} (i) By the standard complex surface theory, when $m>0$,
here the corresponding moduli space $M_\R(h)$ is empty with $h=0$. 
Hence
$SW_\R(\xi)=0$.

(ii) The main issue is that a generic real perturbation $h\in
(\i\Om^2_+)_\R$ may not be generic in $\i\Om^2_+$ (namely the
equivariant transversality fails). However the virtual
neighborhood method can be applied so no generic perturbation is
really necessary to compute $SW(\xi)$. Thus one first uses a
generic real perturbation $h$ to compute $SW_\R(\xi)$.  Then one
applies a suitable neighborhood of the whole moduli space $M(h)$
to compute $SW(\xi)$ (without changing $h$). Furthermore,  when
$m=0$, one can compare the two resulted invariants and prove
$SW_\R(\xi)=SW(\xi)\:\mod\: 2$. The precise argument can be
carried out essentially in the same way
as Ruan-Wang \cite{rw}. \qed\\

Therefore it makes more sense to obtain integer valued real
Seiberg-Witten invariants. The orientability and
orientation of $M_\R(h)$  now come to play, thus
 the line bundle $\det\mr{ind}\de_\R\to\mc{B}^*_\R$ must be invoked.

But first, we have seen that the class $\mu\in H^2(\mc{B}^*, {\bf
Z})$ restricts to a 2-torsion in $H^2(\mc{B}^*_\R, {\bf Z})$,
while $\nu\in H^1(\mc{B}^*_\R, {\bf Z}_2)$ simply does not lift to
the ${\bf Z}$ coefficients. Thus neither class will be useful in
defining integer valued invariants through their pairings with the
possible fundamental class $[M^*_\R(h)]\in H_m(\mc{B}^*_\R,{\bf
Z})$ ($m>0$). In other words, the most likely integer real
Seiberg-Witten invariants come from the virtual dimension $m=0$
real moduli spaces, even for any general almost complex manifold
$X$ admitting real structures. (One might use $H^1(X,{\bf Z})^\si$
to pair the $[M^*_\R(h)]$, but it is not clear how useful the
invariant will be.)

For the rest of the paper, we will  consider all spin$^c$ structures
with virtual dimension 0, unless specifically indicated otherwise.

\begin{thm}\label{prsw}
Fix orientations on  $H^1_\R(X,{\bf R}), H^+_\R(X,{\bf R})$. If
$H^1(X,\R)$ is trivial or more generally if the $\si$
anti-invariant part $H^1(X,\R)^-$ is trivial, then the associated
real Seiberg-Witten invariant is a well-defined integer, possibly
chamberwise when $b^+=1$.
\end{thm}

\nt{\bf Proof}. Consider the usual reducible wall
 $W=\i c^+ +\mr{Im}d^+\subset\i\Om^2_+$,
consisting of all perturbations whose Seiberg-Witten equations
contain reducible solutions.
 It is well-known that the map
 \begin{equation}\label{lk}
\mc{B}^L=\mc{A}/\mc{G}\to W, A\mapsto F^+_A
\end{equation}
is a trivial fibration with fiber the torus $\mc{T}=H^1(X,
\i\R)/H^1(X,2\pi\i{\bf Z})$, see \cite{sa} for example. In particular,
$\mc{B}^L$ is homotopic to $\mc{T}$, since $W$ is contractible.
  Similarly the real version says that
 $\mc{B}^L_\R$ is homotopic to the fixed torus
$\mc{T}_\R=H^1(X,\i\R)^\si/H^1(X,2\pi\i{\bf Z})^\si$. By Lemma 
\ref{hog},
$\dim\mc{T}_\R=\dim H^1_\R(X,\i\R)=\dim H^1(X,\R)^-$.
 Hence $\mc{B}^L_\R$ is contractible under the assumption in the
theorem. Thus the bundle $\det\mr{ind}\Dr_{A,\R}\to \mc{B}^L_\R$
is trivial and oriented. By Theorem \ref{q}, $\det\mr{ind}\de_{\R}$
is trivial and oriented, based on the orientations of
$H^1_\R(X,{\bf R}), H^+_\R(X,{\bf R})$.

Let $M^*_\R(h)\subset\mc{B}^*_\R$ be a regular real Seiberg-Witten
moduli space associated with a generic real perturbation $h\in
(\i\Om^2_+)_\R\backslash W_\R$. By assumption, $\dim M_\R=0$.
Hence at any point $[A,\Phi]\in \dim M_\R$,
$\ker\de_\R=\mr{coker}\de_\R=\{0\}$ and the fiber
$\det\mr{ind}\de_\R$ over $[A,\Phi]$ has a canonical orientation.
As in the proof of Theorem \ref{q}, there is a topological line
bundle $\det\mr{ind}\de_\R$ over the continuous path $[A,t\Phi],
0\leq t\leq 1$ in $\mc{B}_\R$. Any trivialization of the bundle
yields a unique correspondence between  orientations of the fibers
over $[A,0], [A,\Phi]$. Then we define sign$[A,\Phi]=1$ if the
canonical orientation over $[A,\Phi]$ matches the orientation of
the fiber $\det\mr{ind}\de_\R= \det\mr{ind}\Dr_{A,\R}$ over
$[A,0]\in \mc{B}^L_\R $; otherwise define sign$[A,\Phi]=-1$. Then
the real Seiberg-Witten invariant is defined to be the algebraic
sum $\sum\mr{sign}[A,\Phi]$ over all points in $M_\R$. That the
sum is independent of $h$ (on each chamber if $b^+=1$) follows
from the standard cobordism argument, since the sign function
sign$[A,\Phi]$ is continuous and $\mc{B}^L_\R$ is certainly
connected. \qed

Note the theorem still holds for any almost complex manifold $X$
(with $b^+=1$ replaced by $b^+_\R=1$).

To seek an immediate application of the theorem, we consider a real
version of the
Thom conjecture. Let $\Si$ be a closed oriented surface with
a smooth orientation-reversing involution $\tau$. One can show
that the fixed point set  $\Si^\tau$ consists of disjoint circles
and $\Si\backslash\Si^\tau$ has at most two components, see for
example \cite{wil}.
Set $k_\Si$ to be the number of such circles. Call $\tau$ or
$\Si^\tau$ {\em dividing} if $\Si\backslash\Si^\tau$ has exactly
two components.
In this case, let $g^+_{\Si}$ denote the genus of either component.

\begin{cor}\label{tho}
Suppose $\Si\hookrightarrow{\bf C}{\bf P}^2$ is
embedded smoothly and equivariantly with respect to $\tau$ and
the complex conjugation on ${\bf C}{\bf P}^2$. Assume
$\tau$ is dividing and $[\Si]\in H_2({\bf C}{\bf P}^2)$ is 
also represented
by an algebraic curve $C$ of degree $d>2$.
Then $2g^+_{\Si}+k_\Si\geq \frac{(d-1)(d-2)}{2}+1$. In addition, if
$C$ is a dividing real curve in ${\bf C}{\bf P}^2$ and
$k_\Si=k_C$ (the number of ovals in $C_\R=C\cap{\bf R}{\bf P}^2$),
then $g^+_\Si\geq g^+_C$.
\end{cor}

\nt{\bf Proof}. This is essentially an adaptation of
the Kronheimer-Mrowka argument \cite{km} to
our real Seiberg-Witten solutions.

Let $X={\bf CP}^2\# d^2{\ov{\bf CP}}^2$ be a blown-up at
$d^2$ points in $\Si_\R=\Si\cap{\bf RP}^2$ and $\wt{\Si}$
be the internal connected sum with the $d^2$ real
exceptional spheres $E_i$.
Clearly $X$ carries a real structure under which
 $\wt{\Si}$ is invariant. This real structure has a canonical
anti-holomorphic lifting on the line bundle $L=3H-E$, where $H$ is
the hyperplane divisor of ${\bf CP}^2$ and $E=\sum E_i$. Thus the
canonical spin$^c$ structure $S$ on $X$ with determinant $L$
admits a real lifting. By Theorem \ref{prsw} above, the real
Seiberg-Witten invariant of $S$ is well-defined on the two
chambers. The standard argument from Taubes \cite{t1} shows that
the real Seiberg-Witten invariant is 1 on the main chamber, since
the solution from \cite{t1} for a large real perturbation is also
real.

Choose an invariant metric on $\wt{\Si}$ with a constant scalar
curvature $s_0$. Since the real Seiberg-Witten invariant is
non-trivial, by the argument of \cite{km} there is a real
Seiberg-Witten solution $(A,\Phi)$ on $X$, satisfying $|F_A|\leq
-2\pi s_0$ in a neighborhood of $\wt{\Si}\subset X$. Let $\Si^+$
be one component of $\Si\backslash \Si^\tau$; similarly define
$\wt{\Si}^+$. Because $A$ hence $F_A$ is (anti) invariant under
the real structure, we have the following calculations:
$$3d-d^2=c_1(L)[\wt{\Si}]=2\int_{\wt{\Si}^+}\frac{i}{2\pi}F_A,$$
from which we have
$$-(3d-d^2)\leq2\int_{\wt{\Si}^+}\frac{1}{2\pi}|F_A|
\leq 2\int_{\wt{\Si}^+}(-s_0)=2(g^+_{\wt{\Si}}+k_{\wt{\Si}}-2)$$
where the last equation is the Gauss-Bonnet formula
on the surface $\wt{\Si}^+$ with boundary. Note that $\wt{\Si}^+$
is just ${\Si}^+$ connected sum with $d^2$ half disks and
$\pa\wt{\Si}^+$ is $\pa\wt{\Si}^+$ connected sum with $d^2$
semi-circles. Hence
$g^+_{\wt{\Si}}=g^+_{{\Si}}, k_{\wt{\Si}}=k_{\Si}$.
From the computations above,
one arrives at
\begin{equation}\label{gih}
g^+_\Si+k_\Si\geq \frac{(d-1)(d-2)}{2}+1.
\end{equation}

For a real dividing algebraic curve $C\subset{\bf CP}^2$, the Euler
characteristic satisfies $\chi_{_{C}}=2\chi_{_{C^+}}$. In 
terms of genus,
this translates into
$2-2g_{_C}=2(2-2g^+_C-k_C)$, which leads to
\begin{equation}\label{ji}
g^+_C+k_C= g_{_C}+1=\frac{(d-1)(d-2)}{2}+1.
\end{equation}
If $\Si$ is confined  by $k_\Si= k_C$, then the last equation
implies $g^+_\Si\geq g^+_C$
in view of (\ref{gih}).\qed

\nt{\bf Remark.}
(1) In the case that $\tau$ is non-dividing, the corollary remains 
to be true
if one  replaces $g^+_{\Si}$ with the handle number of the quotient
surface $\Si/\tau$, which is a non-orientable surface
with boundary consisting of $k_{\Si}$ circles.

(2) From (\ref{ji}), one has the Harnack inequality
$$k_C\leq \frac{(d-1)(d-2)}{2}+1$$
which gives the upper bound for the number of ovals in any real
algebraic curves $C_\R$ of degree $d$. It seems reasonable to
conjecture the inequality holds true for any smooth equivariantly
embedded surface $\Si\subset{\bf CP}^2$:
$$k_\Si\leq \frac{(d-1)(d-2)}{2}+1,$$
 as long as $[\Si]=[C]=dH$.\\

Without the assumption $H^1(X,\R)=0$,
Theorem \ref{prsw} can be generalized as follows.

\begin{thm}\label{isw}
Fix orientations on  $H^1_\R(X,{\bf R}), H^+_\R(X,{\bf R})$.
If the determinant $L=\det S^+$ of the spin$^c$ structure has its
Chern class $c_1(L)\in H^2(X,{\bf Z})$ divisible by 4,
then the real Seiberg-Witten invariant is a well-defined integer
(chamberwise when $b^+=1$).
\end{thm}

\nt{\bf Proof.} It is enough to show that the line bundle
$\det\mr{ind}\Dr_{A,\R}\to \mc{B}^L_\R$
is orientable with a unique orientation. Then
by Theorem \ref{q}, $\det\mr{ind}\de_{\R}\to\mc{B}^*_\R$ is oriented.
Furthermore, the
real Seiberg-Witten invariant can be constructed  exactly the 
same way as
the proof of Theorem \ref{prsw}.

Fix a base connection $A_0\in\mc{B}^L_\R$, the fiber of the map
(\ref{lk})
over $0\in W$ is naturally diffeomorphic to
$\mc{T}=H^1(X,\i\R)/H^1(X,2\pi\i{\bf Z})$. Thus we have
a complex line bundle $\eta\to \mc{T}$, using the diffeomorphism
to pull back $\det\mr{ind}\Dr_A$. To prove the theorem
we need to show that the associated real
line bundle $\eta_\R\to \mc{T}_\R$ is oriented uniquely.

Since $X$ is K\"{a}hler, $H^1(X,{\bf C})=H^{1,0}\oplus H^{0,1}$.
It follows that $H^1(X,\i\R)$ is naturally isomorphic to $H^{0,1}$
as real vector spaces. This endows a natural complex structure on
$H^1(X,\i\R)$ and hence on $\mc{T}$. Then $\mc{T}$ becomes the
Picard variety of degree zero holomorphic bundles on $X$. Since
$\si$ is a real structure on $X$, its induced map on the complex
torus $\mc{T}$ is now a real structure as well. (Indeed $\si$
induces a real structure on $H^{0,1}$ as seen before.)

Fix a $\si$-compatible complex basis on $H^1(X,\i\R)$ from that 
on $H^{1,0}$.
The tangent bundle of $\mc{T}$ is naturally isomorphic to the trivial
bundle $\mc{T}\times H^1(X,\i\R)$. Hence $\mc{T}$ carries a
natural spin structure that is real compatible with $\si$ in 
the sense of
Definition \ref{tsi}. In turn this spin structure will determine
a canonical square root of $\eta$ if $\eta$ has one. It is a classical
fact that square roots of $\eta$ are in one-to-one correspondence
with spin structures on $\eta$, for example from \cite{a}.
Hence, assuming $\eta$ has a square root, there is a well-defined
spin structure on $\eta$, which is real compatible with $\si$, because
the spin structure on $\mc{T}$ is so. Applying the main result
in Wang \cite{w3}, we see that the real line bundle $\eta_\R$ is
orientable with a well-defined orientation.

It remains to show that $\eta$ has a square root, namely
$c_1(\eta)\in H^2(\mc{T},{\bf Z})$ is divisible by $2$. Apply the
 Atiyah-Singer family index theorem to the universal
spin$^c$ bundle $\mb{S}$ on $\mc{T}\times X\subset\mc{B}^L\times X$
from Subsection \ref{ru}. Thus
  $ch(\mr{ind}\Dr_A)=\int_X\hat{A}(X)\mr{exp}(\mb{L}/2)$.
As in \cite{ll, ot}, one computes the integral routinely, obtaining
\begin{equation}\label{si}
c_1(\mr{ind}\Dr_A)=\frac{1}{2}\sum_{i<j}<c_1(L)\al_i\al_j,
[X]>\beta_i\beta_j
\end{equation}
where $\{\al_i\}$ is any basis of $H^1(X,{\bf Z})$ and
$\{\beta_i\}$ is the induced dual basis in $H^1(\mc{T},{\bf Z})$.
From our assumption $4 | c_1(L)$, $c_1(\mr{ind}\Dr_A)$ is then
divisible by $2$, so is $c_1(\eta)$ and the proof is finished. In
fact let us take any complex line bundle $K$ on $X$ with $K^2=L$.
From $w_2(X)\equiv c_1(L)\equiv0$ mod $2$, $X$ is spin. Since
$2|c_1(K)\in H^2(X,{\bf Z})$, $c_1(K)$ is a characteristic
element. Thus there is a spin$^c$ structure on $X$ with
determinant $K$. Repeating the above argument for this new
spin$^c$ structure, one sees  the analogy of formula (\ref{si})
implies that $\frac{1}{2}\sum_{i<j}<c_1(K)\al_i\al_j,[X]>$ are all
integers for any $i<j$. It follows that $c_1(\mr{ind}\Dr_A)$ is an
even class. \qed

Note that for the theorem, it is not enough to assume only $2| c_1(L)$,
because then the bundle $K$ in the last part of the proof will not be
characteristic and $c_1(\mr{ind}\Dr_A)$ may not be divisible by 2.

Theorem \ref{isw} can be extended to symplectic manifolds such that
$b_1(X)$ is even.

 To give some examples with the Chern class
$c_1(X)$ divisible by $4$, we can take a product of Riemann
surfaces, $\Si_g\times \Si_h$, with odd genera $g, h$. Here both
$\Si_g$ and $\Si_h$ carry real structures. If $c_1(X)$ is
divisible by $4$, we can get additional examples by taking any
branched cover of $X$ along a branched locus $C\subset X$ that is
preserved by the real structure and such that $4|[C]$.

\section{Seiberg-Witten projection maps}\label{yup}

In the initial part of the section we work with the most general
real set-up, assuming only that $(X,\si)$ is any smooth 4-manifold
with involution and  $P_{sp}$ is a spin$^c$ structure that is
endowed  with a real compatible lifting of $\si$ in the sense of
Definition \ref{tsi}. Then the Seiberg-Witten equations inherit a
real structure as in Proposition \ref{rsw}. So far we have studied
the real Seiberg-Witten moduli spaces by studying the ambient
configuration space $\mc{B}^*_\R$, with the approach
 parallel to the standard theory. In this section, we will 
shift our focus
 and investigate the moduli spaces directly without going over 
$\mc{B}^*_\R$.
 More precisely let $Q=\i\Om^2_+\backslash W,
Q_\R=(\i\Om^2_+)_\R\backslash W_\R$ denote
 the complements of the (real) reducible walls. Then
we will analyze systematically the Seiberg-Witten projection and
its real version:
$$\pi:{\bf M}\to Q, \pi_\R:{\bf M}_\R\to Q_\R,$$
where ${\bf M}, {\bf M}_\R$ are the parameterized (irreducible)
moduli spaces. Note that both projections are proper smooth
Fredholm maps by the usual compactness theorem. (In comparison,
the full projection ${\bf M}\to\i\Om^2_+$ is only a continuous
proper map, while the restriction to the irreducible ones ${\bf
M}^*\to\i\Om^2_+$ is smooth but not proper.)

\subsection{The structures of critical points and critical values}
In this subsection, we can actually consider an arbitrary
Fredholm index of $\pi$, i.e. the virtual dimension $\mr{ind}\de$
of the moduli space is any integer. In fact, a point of our
approach is to extract possibly additional information from $\pi$
or $\pi_{\R}$ in the case of a negative virtual index where the
usual Seiberg-Witten invariant fails to yield any information.
Compare with Shevchishin \cite{s} where the moduli space of
pseudo-holomorphic curves was studied.

First we consider the general situation.
Let $\mc{E}, \mc{F}$ be Banach bundles over $M$, and
$\ell:\mc{E}\to\mc{F}$ be a Fredholm bundle homomorphism of constant
index $m=\mr{ind}\ell_x, x\in M$. Then
using connections on $\mc{E}, \mc{F}$, one can define a 
pointwise linear map
$$\nabla\ell:T_xM\lra Hom(\ker\ell_x,\mr{coker}\ell_x)$$
for each $x\in M$, which is actually independent of the
connections chosen, see Lemma 1.3.1 of Shevchishin \cite{s}. The
following basic result is used on page 50 of \cite{s}  without
proof:

\begin{lem}\label{shev}  Let $C(l)=\{x\in M\mid\dim
\mr{coker}\ell_x=l\}$ for a fixed integer $l\geq0$. If $\nabla
\ell_x$ is surjective for all $x\in C(l)$, then $C(l)\subset M$ is
a submanifold of codimension $(m+l)l$.
 \end{lem}

 \nt{\bf Proof.} We sketch for the case where
$\mc{E}=M\times U, \mc{F}=M\times V$ are trivial product bundles,
which is what we require in our applications. The general case can
be dealt with using suitable modifications.

Consider the Banach space $\mr{Fred}(U,V)_m$ of all Fredholm
operators of index $m$. The subset
$$W=\{g\in\mr{Fred}(U,V)_m\mid \dim\mr{coker} g=l\}$$
is a submanifold of codimension $(m+l)l$. The map $\ell$ becomes
$M\to\mr{Fred}(U,V)_m$  and $\nabla\ell=p\circ d\ell\circ i$,
where $i:\ker\ell_x\hookrightarrow U, p:V\to\mr{coker}\ell_x$. One
may check that $\ell$ is transversal to $W$ iff $\nabla\ell$ is
surjective on $N(l)$, by noting that the tangent space of $W$ is
$$T_gW=\{h\in Hom(U,V)\mid h \text{ maps } \ker g 
\text{ to  im}g\}.$$
It follows then from the usual transversality theorem that
$C(l)=\ell^{-1}(W)$ is a submanifold of codimension $(m+l)l$.\qed

\begin{cor} If $\nabla\ell$ is always surjective at any point
$x\in M$, then $M$ is stratified by submanifolds $C(l), l=0, 1,
\cdots$.
 \end{cor}

Return to our parameterized Seiberg-Witten moduli space 
${\bf M}$ and the
projection into the perturbation space
 $\pi:{\bf M}\lra Q$ as in Subsection \ref{prel}.
  Let ${\bf C}$
 denote the critical point set of $\pi$ and
 ${\bf C}(l)=\{{\bf x}\in{\bf M}\mid \dim\mr{coker}D\pi_{\bf x}=l\}$.

\begin{thm}\label{cr} For each $l=0, 1,2,\cdots$,
${\bf C}(l)\subset {\bf M}$ is a Banach submanifold of codimension
$kl$, where $k=\mr{ind}D\pi+l$.
 \end{thm}

\nt{\bf Proof.} From Corollary (\ref{dee}), it is the same to show
that ${\bf C}(l)=\{{\bf x}\in {\bf M}\mid \dim\mr{coker}\de_{\bf
x}=l\}$ is a codimension $kl$ submanifold of ${\bf M}$. We can of
course view
$$\de:{\bf M}\lra \mr{Fred}(U, V),$$
where $U=\i\Om^1\oplus\Ga(S^+),
V=\i\Om^0\oplus\i\Om^2_+\oplus\Ga(S^-)$ (The suitable Sobolev
spaces are suppressed without harm). To apply Lemma (\ref{shev}),
we need  to show $\nabla\de_{\bf x}$ is surjective.

Let us compute the differential $d\de_{\bf x}:\mr{T}_{\bf x}{\bf
M}\to Hom(U,V)$. Take a point ${\bf x}=(A,\Phi,h)\in{\bf C}(l)$, a
tangent vector $\xi=(a,\phi,k)\in\mr{T}_{\bf x}{\bf M}$ and
$(a',\phi')\in U$. Then we have in $V$ that:
\begin{equation}\label{d3}
 d\de_{\bf x}(\xi)(a',\phi')=
(\i<\phi,\phi'>,Dq_{\phi}(\phi'),2^{-1}a'\cdot\phi+2^{-1}a\cdot\phi').
\end{equation}
  Consider $\nabla\de_{\bf x}:\mr{T}_{\bf x}{\bf M}\to 
Hom(\ker\de_{\bf
x}, \mr{coker}\de_{\bf x})$, with $\nabla\de_{\bf x}(\xi)$ equal
to the composition
 \begin{equation}\label{xx}
 \begin{array}{ccccc}
 &d\de_{\bf x}(\xi)&&p&\\
 \ker\de_{\bf x}\hookrightarrow U&\to&V&\to&\mr{coker}\de_{\bf x}.
 \end{array}\end{equation}

 We need to show that by choosing $\xi$ suitably, $\nabla\de_{\bf
 x}(\xi)$ can realize all linear functions
 $f(a',\phi')$ from $\ker\de_{\bf x}$ to $\mr{coker}\de_{\bf x}$.
 Note that each of the three components of $d\de_{\bf x}$ from
 (\ref{d3}) is non-degenerate bilinear in the two sets of
 variables
 $\{a,\phi\}$ and $\{a',\phi'\}$. Hence each component can realize
 all linear functions of one set of variables $\{a',\phi'\}$ when
 the other set $\{a,\phi\}$ is suitably chosen. Of course this
 does not mean that all three components can simultaneously
 realize arbitrarily given three functions. However, after 
composing with
 the
projection map $p$, only two components are actually independent.
Moreover, when we restrict to $\ker\de_{\bf x}$, the two variables
$a', \phi'$ are not independent either. Therefore, essentially
just one independent variable from the set $\{a, \phi\}$ is needed
in order for the composition (\ref{xx}) to realize all linear
functions $f$ as indicated above. It would seem that we have a
redundant variable from $\{a,\phi\}$, but remember
$\xi=(a,\phi,k)\in \mr{T}_{\bf x}{\bf M}$ must satisfy two
equations
$$D^2_{\bf x}(a,\phi)=(k,0), (D^1_{\bf x})^*(a,\phi)=0$$
according to Proposition (\ref{dpi}). So actually we only have one
essentially independent variable available from $\xi$, and this is
good enough here.\qed\\

\nt{\bf Remark.} Even when $X$ is a K\"{a}hler manifold, $\ker
D\pi, \mr{coker}D\pi$ may be of odd dimensions at a non-trivial
perturbation $h$.

Next we take up the set up with a real structure,
so we have the real Fredholm map $\pi_{\R}:{\bf M}_{\R}\to Q_{\R}$.
Let ${\bf C}_{\R}$ be the critical point set of $\pi_{\R}$ and
${\bf C}_{\R}(l)$ be the subset of points at
which $\dim\mr{coker}D\pi_{\R}=l$. Thus ${\bf C}_{\R}(0)$ is the set
of regular points of $\pi_\R$.
The real version of Theorem \ref{cr} holds under the same proof:

\begin{thm}\label{crr}For each $l$, 
${\bf C}_{\R}(l)\subset{\bf M}_{\R}$
is a Banach submanifold of co-dimension $l(\mr{ind}D\pi_{\R}+l)$.
\end{thm}

In particular, when the virtual dimension $\mr{ind}D\pi_{\R}=0$,
the subset ${\bf C}_{\R}(1)$ is a co-dimension 1 submanifold in
${\bf M}_{\R}$.

\subsection{Degree of Seiberg-Witten projection map}
In this subsection, we study the projection $\pi_\R$ from a
functional analytic point of view. Suppose in general that $f:M\to
N$ is a proper smooth Fredholm map of index $0$ between two Banach
manifolds. In order to define an integer degree of $f$, the most
natural approach is to impose certain oriented manifold structures
on $M, N$ and require $f$ to preserve these structures. The only
subtlety here is that the general linear group $GL(E)$ of an
$\infty$-dimensional Hilbert space $E$ is contractible, thus
connected, by a classical result of Kuiper. Hence, one needs to
reduce the structure group of $TM,TN$ to the smaller subgroup
$GL_c(E)$ of compact linear isomorphisms which has two connected
components, so that the orientability may be imposed. This was the
approach initiated by K.D. Elworthy and A.J. Tromba in the 1970s.

More recently, Fitzpatrick, Pejsachowicz, and Rabier \cite{fpr}
realized that the orientability of $M,N$ is often un-natural to
impose and not necessary either for the sole purpose of defining a
degree for $f$. Instead, all needed is the orientability of the
map $f$ itself. In \cite{fpr}, they introduced the parity  of $f$
along a path  with two ends at regular points of $f$. This is a
functional analytic concept which involves parametrices and the
Leray-Schauder mod-2 degree. Then $f$ is called {\em orientable}
if the parity is always $1$ along any loop.

On the other hand,  the geometric point of view is to characterize
 the orientability of $f$ as that of the determinant line bundle
$$\det f=\wedge^{\mr{max}}\ker Df\otimes(\wedge^{\mr{max}}
\mr{coker}Df)^*$$
over $M$ using the Fr\'echet derivative $Df:TM\to f^*TN$.
It is proved in \cite{w2} that the two kinds of 
orientability mentioned
above are actually equivalent. Namely, $\det f$ is a 
trivial line bundle
iff $f$ is orientable in the  sense of \cite{fpr}.
 Let $C_f\subset M$ denote the set of critical points
where
$\mr{coker} Df$ is 1-dimensional
and $R_f$ the set of regular points of $f$. Then the
equivalence in turn leads
to the following (see \cite{w2}):

\begin{pro}\label{wor}
Suppose that $C_f$ is a co-dimension 1 submanifold of $M$ and
$R_f\not=\emptyset$. Then the line bundle $\det f$ is trivial  iff
there is a continuous sign function $\epsilon: R_f\to\{\pm1\}$,
such that for any path $\gamma\subset M$ with both ends in $R_f$
and transversal to $C_f$, the sign $\epsilon$ will change whenever
$\gamma$ crosses $C_f$.
 \end{pro}

Naturally the parity of $f$ along a path between two regular points
can now be determined by $\epsilon$.
Each  $\epsilon$ is called an {\em orientation} of $f$ in \cite{fpr}.
By \cite{w2}, this corresponds canonically to an orientation 
of $\det f$.
 Proposition \ref{wor} gives a convenient criterion for 
the orientability
and orientation of $\det f$ in terms of  signs at regular points only.

From here on we understand that $f$ is {\em oriented} if $f$ carries
 a sign function $\epsilon$ as in Proposition \ref{wor}.
 Then the integer degree is defined to be
$$\deg f=\sum_{x\in f^{-1}(y)}\epsilon(x),$$
where $y\in N$ is a regular value.

Recall from \cite{fpr} that an oriented homotopy is a smooth
Fredholm map $H:M\times[0,1]\to N$ that carries an orientation.
Using determinant bundles, it is easy to see that a homotopy $H$
is orientable (oriented) iff some section $H_t:M\times\{t\}\to N$
is orientable (oriented respectively). Note that $\det f$ is not
exactly homotopy invariant in the usual sense; instead we  should
utilize the following (see \cite{fpr}):

\begin{pro}\label{depro} Suppose $f$ is an oriented Fredholm
map of index zero.

\nt{\em (Homotopy Invariance)}  The degree
$\deg f$ is invariant under any proper and oriented homotopy $H$.
Hence the absolute
value $|\deg f|$ is homotopic invariant regardless of orientation.

\nt{\em (Reduction)} If $P\subset N$ is a submanifold transversal
to $f$, then the restriction $f|_P: f^{-1}(P)\to P$ is a Fredholm
map with an induced orientation. Moreover, $R_f\cap f^{-1}(P)$ gives
all regular points of $f|_P$ and consequently $\deg f=\deg f|_P$.
\end{pro}

We now return to our Seiberg-Witten projection
$\pi_{\R}:{\bf M}_{\R}\to Q_{\R}$, assuming the virtual dimension
is zero.
The main point is that
$\pi_{\R}:{\bf M}_{\R}\to Q_{\R}$ can be orientable,
although $\det\mr{ind}\de_\R\to \mc{B}_\R$
 may well be  non-trivial, making Section \ref{rss}
inapplicable. (This is in analogy  with \cite{we} where
only rational curves are given suitable signs). In other words
we can expand the definition from Section \ref{rss}:

\begin{df}\label{ssr}
When $\pi_{\R}$ is oriented, the real Seiberg-Witten invariant
$SW_\R(P_{sp})$
is defined to be the degree of $\pi_{\R}$.
\end{df}

By Proposition \ref{depro}, with fixed orientations on $H^1_\R(X,\R)$
and $H^+_\R(X,\R)$, $SW_\R(P_{sp})$ is independent of metrics on $X$.
Without fixing the orientations, the absolute value $|SW_\R(P_{sp})|$
is still well-defined.

To detect the orientability, from Theorem \ref{crr} and
Proposition \ref{wor}, it is enough to give a continuous sign
assignment $\epsilon$ at the regular points of $\pi_\R$ such that
$\epsilon$ changes whenever crossing the submanifold ${\bf
C}_\R(1)$. In general it is still rather difficult to find a
suitable $\epsilon$. Nonetheless one immediate result within reach
is a real blow up formula, which we describe next. Let
$\hat{X}=X\#\ov{\bf CP}^2$ be the blow-up of $X$ at a real point.
Then $\si$ extends smoothly over $\hat{X}$ as an involution, which
further lifts to the spin$^c$ bundle $\hat{P}_{sp}$ on $\hat{X}$.
 Let $\hat{\pi}_\R:\hat{{\bf M}}_\R\to\hat{Q}_\R$ be the 
real Seiberg-Witten
projection in the spin$^c$ structure on $\hat{X}$. Here is the 
real version
of the usual blow up formula, the counter part of which is much 
harder to
prove for real rational curves in \cite{we}.
\begin{thm}\label{bl}
If $\pi_\R:{{\bf M}}_\R\to{Q}_\R$ is an orientable Fredholm map,
then so is $\hat{\pi}_\R:\hat{{\bf M}}_\R\to\hat{Q}_\R$. Moreover
an orientation of $\pi_\R$ induces one for $\hat{\pi}_\R$ and
the real Seiberg-Witten invariant remains the same:
$SW_\R(P_{sp})=SW_\R(\hat{P}_{sp})$.
\end{thm}

\nt{\bf Proof.} One just needs to make sure that the usual proof
can be carried out equivariantly with respect to our real
structures. Let $S^2={\bf CP}^1$ be given the standard complex
conjugation. The degree $-1$ line bundle on $S^2$ has a natural
real lifting, which preserves the standard Hermitian fiber metric.
Thus the disk bundle $N$ inherits the real structure, which of
course is the restriction of the complex conjugation to the
neighborhood of ${\bf CP}^1\subset\ov{\bf CP}^2$. Attach a long
cylinder $[1, r]\times S^3$ to the boundary $\pa{N}=S^3$ and let
$N_r$ denote the resulted manifold with the extended real
structure.

Fix a small 4-disk $D\subset X$ at the real blow up point in 
$X$ and let
$D_r$ denote manifold with a long cylinder attached. Attach 
this cylinder
 as well as the infinite cylinder $[1, \infty)\times S^3$ to 
the boundary
$S^3 =\pa(X\backslash D)$ so we get two more manifolds 
$X_r, X_\infty$. Note
that $D_r, X_r, X_\infty$ all inherit real structures from $X$.

Fix a large enough $r$ and diffeomorphisms $X\approx X_r\cup D_r,
\hat{X}\approx X_r\cup N_r$ (glue the long necks together).
Without loss of generality we assume all perturbations on $X$ and
$\hat{X}$ have compact supports on $X_r$, namely they are trivial
near the blow up point and the exceptional curve $S^2$
respectively. Thus we have identified the (real) perturbation
spaces $Q_\R\approx \hat{Q}_\R$. Via the  standard gluing process,
every irreducible Seiberg-Witten solution on $X$ and on $\hat{X}$
both correspond to a unique finite energy solution on $X_\infty$.
Thus we have the usual diffeomorphism ${\bf M}\approx\hat{{\bf
M}}$ between the parameterized moduli spaces. The gluing between
$X_r$ and $N_r$ requires the use of a canonically defined
reducible solution $(A_0,0)$ in the spin$^c$ structure over $N$.
One checks easily that this solution is equivariant with respect
to the real structure. Thus one has a diffeomorphism ${\bf
M}_\R\approx \hat{\bf M}_\R$ by restriction. Since $\pi_\R,
\hat{\pi}_R$  commute under the previous diffeomorphisms ${\bf
M}_\R\approx \hat{\bf M}_\R, Q_\R\approx \hat{Q}_\R$, the
orientability of one certainly implies that of the other. \qed

In particular if the spin$^c$ structure $P_{sp}$ satisfies the
condition in Theorem \ref{isw}, then the blow up spin$^c$
structure $\hat{P}_{sp}$ on ${X}\#\ov{{\bf CP}}^2$ has an integer
real Seiberg-Witten invariant. Note the determinant of
$\hat{P}_{sp}$ is no longer divisible by $4$; thus Theorem
\ref{isw} is not applicable. Some of the standard applications of
the blow up formula can be readily extended to our real case.

 Another observation to make is about reductions. Assume here
that $(X,\om)$
is a symplectic real 4-manifold.
 Recall $Q_\R\subset(\i\Om^0)_\R\cdot\om\oplus(\i\Om^{0,2})_\R$.
Let $\pi^1_\R,\pi^2_\R$ be
respectively the compositions of $\pi_\R$ with the projections onto
the two factors. It is straight forward to show the following

\begin{pro}\label{tran}
The map $\pi_{\R}$ is orientable (oriented)
iff $\pi_{\R}^1, \pi_{\R}^2$ are
orientable (oriented respectively).
\end{pro}

The Fredholm map $\pi_{\R}^1:{\bf M}_\R\to(\i\Om^0)_\R\cdot\om$ 
involves
the generalized Taubes perturbation \cite{t1}, whereas 
$\pi^2_{\R}:{\bf
M}_\R\to(\i\Om^{0,2})_\R$ involves a Witten type 
perturbation \cite{wi}.
With either perturbation, the Seiberg-Witten equations can be
decomposed nicely, and the orientability can be sorted out in
special situations.

\subsection{The non-orientable case: chamberwise invariants}
The set-up here is an almost complex 4-manifold $(X, J)$  
together
with a real structure. It is a basic fact that $J$ maps
isomorphically the tangent space $TX_\R$ to the normal space of
$X_\R$. Applying the same to the domain and range of
$\de=\de_{A,\Phi}$ and in view of Corollary \ref{rdee}, we see
easily that regular points of $\pi_\R$ are
 real regular points of $\pi$. Thus
 we can restrict the usual orientation of $\pi$
to $\pi_{\R}$. Namely the orientation of the standard
Seiberg-Witten theory gives the continuous sign map $\ep_c$
for $\pi$, and by restriction we obtain a sign map $\ep$ defined
at the regular points of $\pi_{\R}$.
(Note that $D\pi$ is not complex linear, hence not all signs of
 $\ep_c$ are positive.)
Then one can apply the criterion in Proposition \ref{wor} together
with Theorem \ref{crr} and seek to determine when $\ep$ is in
fact an orientation for $\pi_{\R}$. We believe this
should work for a class of real almost complex 4-manifolds 
that include
cases in Theorems \ref{prsw} and \ref{isw}, although we have not
checked  the details. (The last claim is  essentially in view
of the deformation of $\de_\R$ to the linear operator $\de^L_\R$.)

What we are interested more is about the opposite case that the
sign map $\ep$ is {\em not} an orientation for $\pi_\R$, as it  will
bring up  new  geometry to study. Specifically
let $Z, T$ be respectively the sets of regular
values and critical values of $\pi_{\R}$. Since $\pi_{\R}$ is
proper, $Z$ is open and dense in $Q_{\R}$ by the Sard-Smale
theorem. Call  connected components of $Z$ the {\em chambers},
which are divided by the {\em wall} $T$.

Take any regular value $h\in Z$, we can count
the signed points in $\pi^{-1}_{\R}(h)$ using our map $\ep$.
Note that $\pi^{-1}_\R(Z)$ is generally a proper subset of the
regular point set ${\bf C}_\R(0)$ of  $\pi_\R$, so we could require
the map $\ep$ be defined in a smaller set than ${\bf C}_\R(0)$.
Obviously the resulted number is independent of regular values in
the same chamber.
Hence it makes sense to define the chamberwise real Seiberg-Witten
invariant for a real almost complex 4-manifold.  For example, in
 the Taubes chamber that contains $\i r\om, r$ a large constant,
 the real Seiberg-Witten invariant takes value $\pm1$, since the only
(regular) Seiberg-Witten solution from Taubes' argument \cite{t1}
is also real.

In the non-orientable $\ep$ case, the real Seiberg-Witten
invariant will vary from chamber to chamber.
The pattern and distribution of the invariant then become the
new geometry to investigate. The essential issue is to give
a ``wall crossing formula'' that describes
 the change between two neighboring chambers in $Z$.
More precisely take any path
$\Gamma=\{\ga(t)\in Q_\R, -1\leq t\leq 1\}$ that is transversal to
$\pi_\R$,
so that $\Gamma'=\pi^{-1}_\R(\Gamma)$ is a submanifold consisting
of finitely many arcs. Suppose all points, except $\ga(0)$,
in $\Gamma$ are regular values and $\ga(-1), \ga(1)$ belong to
different chambers.
We need to examine the restriction $\ti{\pi}:
\Gamma'\to \Gamma$ of $\pi_\R$. Since $\pi_\R$ is
transversal to the 1-dimensional $\Gamma$, $\dim\mr{coker}D\pi_\R$
is at most 1
at any point in $\Gamma'$. Hence at any
critical point of $\ti{\pi}$, $\dim\mr{coker}D\pi_\R$ is 
exactly 1. The
converse is also true; therefore the critical point set of
$\ti{\pi}$ equals ${\bf C}_\R(1)\cap\Gamma'.$

Set $h_{\pm}=\ga(\pm1), q=\ga(0)$. The pre-image points
$\ti{\pi}^{-1}(h_{\pm})$ all carry
signs according to $\ep$. To describe the invariant change between
the two chambers means to compare the two sets of signs here.
Take an arc component $\eta$  of $\Gamma'$.
Along $\eta$, the only possible critical point of $\ti{\pi}$
is $p\in\eta\cap \ti{\pi}^{-1}(q)$. If $p$ is not a critical
point, of course
the two ends of $\eta$ should have the same sign by continuity of
$\ep$. Otherwise we can determine its type:

\begin{pro}\label{cham}
If $p$ is a critical point of $\ti{\pi}$ along $\eta$, then $p$ is
a non-saddle point. Namely $\ti{\pi}$ has either a local maximum
or a local minimum at $p$, under suitable re-parameterizations of
$\eta$ and $\Gamma$.
\end{pro}

\nt{\bf Proof.} Here we adapt a Kuranishi type argument of a
finite dimensional reduction (which  also reflects how the
Leray-Schauder mod-2 degree is defined). In essence, this is due
to the fact that the only non-linear part of the Seiberg-Witten
equations (\ref{SW}) is the quadratic term $q(\Phi)$.

Let $L=D\pi_\R(p): T_p{\bf M}_\R\to T_q Q_\R$ be the differential
at $p$.
Up to diffeomorphisms and locally around $p\in {\bf M}_\R$, we can
decompose
\begin{equation}\label{lst}
\pi_\R(u,v)=(Lu, \psi(u,v))\in \mr{im}L\times\mr{coker}L,
\end{equation}
 where $(u,v)\in \ker{L}^\perp\times\ker{L}$ and $\psi$ is a function
with $D\psi(0,0)=0$. Recall the critical point $p\in{\bf
C}_\R(1)\cap\Gamma'$, meaning that $\mr{coker}L$ and hence $\ker
L$ are both 1-dimensional spaces. Certainly $\psi$ depends on the
various choices made. But as in the original Donaldson theory, the
quadratic part of the restriction $f(v)=\psi(0,v):\ker L\to
\mr{coker}L$ is intrinsic, namely after re-parameterizations $\ker
L\approx {\bf R}, \mr{coker}L\approx{\bf R}$, we always have
\begin{equation}\label{lt}
f(v)=\pm v^2+O(v^2).
\end{equation}
On the other hand, our spaces $\eta, \Gamma$ are also 1-dimensional.
Applying the Implicit Function Theorem if necessary
and in  view of (\ref{lst}),
we can assume  that locally $\eta=\ker L,
\Gamma= \mr{coker}L$ and $\ti{\pi}=f$. By (\ref{lt}),
$\ti{\pi}$ has a local maximum or minimum at $p=0$.
\qed

It follows from Proposition \ref{cham} that $\ep$ is an orientation iff
 every such an $\eta$ must have opposite signs at its two ends.
Other than Proposition \ref{cham}, we have not yet 
determined any precise
wall crossing formula but conjecture
that the invariant change should be independent of neighboring chambers.

\vspace{6mm}

\nt{\bf Appendix: Real classes and classes of real points}
\vspace{2mm}

We lay down the following useful algebraic set up once for all.
 It has scattered widely in the literature that deals with
real structures.

Assume that $C$ is a set and $\si:C\to C$ an involution.
Write $\ov{x}=\si(x)$ for convenience, where $x\in C$.
Analogously for a group $G$, let $\si':G\to G, g\mapsto\ov{g}$ be
an involution such that
$$\ov{gh}=\ov{g}\ov{h}, \ov{1}=1, \mbox{ for } g, h\in G.$$
(Namely $\si'$ is a group homomorphism.)
Suppose $G$ acts {\em freely} on $C$ and the involution actions are
compatible in the sense that
$$\ov{gx}=\ov{g}\:\ov{x}\; \mbox{ for } g\in G, x\in C.$$
From this, $\si$ and $\si'$ induce an involution
$\si_*$ on the quotient set $B=C/G$.

We need to introduce additional sets.
If $\si,\si', \si_*$ are viewed as real structures, then the set of
{\em real classes} should be $B^{\si_*}:=
\mbox{Fix}(\si_*:B\to B)\subset B$,
while the set of {\em classes of  real points} should be the quotient
$$B_{\R}:=C_{\R}/G_{\R}=\mbox{Fix}\:\si/\mbox{Fix}\:\si'.$$
There is a natural inclusion $B_\R\hookrightarrow B^{\si_*}$.
The main purpose
of the Appendix is to generalize the set $B_\R$ as well as the
inclusion.

Define a subgroup $U=\{g\in G:\ov{g}g=1\}$ of $G$ and its quotient
$\wt{U}=U/\sim$, where $g\sim\ov{h}gh^{-1}$ for some $h\in G$.
Any $g\in U$ yields involutions $\si_g:C\to C, x\mapsto \ov{gx}$ and
$\si'_g: G\to G, h\mapsto g^{-1}\ov{h}g$ which are compatible in the
above sense. (Note that all elements in $\wt{U}$ have order 2.)
One can view $\si_g, \si'_g$ as shifted real structures by $g$.
With these new real structures,  we introduce the set
$$B_g=\mbox{Fix}\:\si_g/\mbox{Fix}\:\si'_g,$$
generalizing that $B_1=B_{\R}$.
Given a class $\xi\in \wt{U}$, we introduce a subset of $B$:
$$B^\xi=\{[x]\in B: \ov{x}=gx \mbox{ for some } g\in \xi\}.$$
With the right set up at hands,  one can verify easily the following
statements.
\vspace{2mm}

\nt{\bf Proposition.} {\em (i) The subsets $B^\xi$, with
$ \xi\in \wt{U}$, are mutually disjoint.

(ii) Clearly each $B^\xi\subset B^{\si_*}$; moreover there is a
natural decomposition:
$$B^{\si_*}=\coprod_{\xi\in\wt{U}}B^\xi=B^{[1]}\coprod
(\coprod_{\xi\not=[1]}B^\xi).$$

(iii) There is a natural bijection $B_{\R}\to B^{[1]}, [x]_{\R}\mapsto
[x]$. In particular, we have an inclusion
$B_\R\hookrightarrow B^{\si_*}$.
More generally, we have a natural bijection $B_g\to B^{[g]},
[x]_{\si_g}\mapsto [x]$, where $g\in U, [g]\in\wt{U}$. Thus we 
can rephrase
the previous decomposition as
$$B^{\si_*}=\coprod_{[g]\in\wt{U}}B_g=B_{\R}\coprod
(\coprod_{[g]\not=[1]}B_g).$$}
\vspace{2mm}

In topological applications,
 one usually expects that $\coprod_{[g]\not=[1]}B_g$ 
constitutes a small
subset of $B^{\si_*}$ relative to $B_\R$.

The proportion has been applied to the real and fixed 
configuration spaces
$\mc{B}^*_\R, \mc{B}^{*\si}$ in Section
\ref{dlin}, in which $\mc{G}$ acts freely on $\mc{C}^*$.

\end{document}